%
%
%
%
\hsize=5in
\baselineskip=12pt
\vsize=19.5cm
\parindent=.5cm
\predisplaypenalty=0
\hfuzz=2pt
\frenchspacing
%
%
\input amssym.def
\def\titlefonts{\baselineskip=1.44\baselineskip
	\font\titlef=cmbx12
	\titlef
	}
\font\ninerm=cmr9
\font\ninebf=cmbx9
\font\ninei=cmmi9
\skewchar\ninei='177
\font\ninesy=cmsy9
\skewchar\ninesy='60
\font\nineit=cmti9
\def\reffonts{\baselineskip=0.9\baselineskip
	\textfont0=\ninerm
	\def\rm{\fam0\ninerm}%
	\textfont1=\ninei
	\textfont2=\ninesy
	\textfont\bffam=\ninebf
	\def\bf{\fam\bffam\ninebf}%
	\def\it{\nineit}%
	}
%
%
\def\frontmatter{\vbox{}\vskip1cm\bgroup
	\leftskip=0pt plus1fil\rightskip=0pt plus1fil
	\parindent=0pt
	\parfillskip=0pt
	\pretolerance=10000
	}
\def\endfrontmatter{\egroup\bigskip}
\def\title#1{{\titlefonts#1\par}}
\def\author#1{\bigskip#1\par}
\def\address#1{\bigskip{\reffonts\it#1}}
\def\email#1{\bigskip{\reffonts{\it E-mail: }\rm#1}}
\def\thanks#1{\footnote{}{\reffonts\rm\noindent#1\hfil}}
\def\abstract{\bgroup\reffonts
	\noindent{\bf Abstract. }}
\def\endabstract{\par\egroup}
\def\section#1\par{\ifdim\lastskip<\bigskipamount\removelastskip\fi
	\penalty-250\bigskip
	\vbox{\leftskip=0pt plus1fil\rightskip=0pt plus1fil
	\parindent=0pt
	\parfillskip=0pt
	\pretolerance=10000{\bf#1}}\nobreak\medskip
	}
\def\proclaim#1. {\medbreak\bgroup{\noindent\bf#1.}\ \it}
\def\endproclaim{\egroup
	\ifdim\lastskip<\medskipamount\removelastskip\medskip\fi}
\def\item#1 #2\par{\ifdim\lastskip<\smallskipamount\removelastskip\smallskip\fi
	{\rm#1}\ #2\par\smallskip}
\def\jtem#1 #2\par{\ifdim\lastskip<\medskipamount\removelastskip\medskip\fi
	{\setbox0\hbox{(F2)\ }\leftskip=\wd0
	\noindent\hskip-\leftskip
	\hbox to\leftskip{\rm#1\ \hfil}\it#2\par}\medskip}
\def\Proof#1. {\ifdim\lastskip<\medskipamount\removelastskip\medskip\fi
	{\noindent\it Proof.}\ }
\def\endproof{\quad\hbox{$\square$}\medskip}
\def\Remark. {\ifdim\lastskip<\medskipamount\removelastskip\medskip\fi
	{\noindent\bf Remark.}\quad}
\def\Remarks. {\ifdim\lastskip<\medskipamount\removelastskip\medskip\fi
	{\noindent\bf Remarks.}\quad}
\def\endremark{\medskip}
\def\Example. {\ifdim\lastskip<\medskipamount\removelastskip\medskip\fi
	{\noindent\bf Example.}\quad}

%
%
\newcount\citation
\newtoks\citetoks
\def\citedef#1\endcitedef{\citetoks={#1\endcitedef}}
\def\endcitedef#1\endcitedef{}
\def\citenum#1{\citation=0\def\curcite{#1}%
	\expandafter\checkendcite\the\citetoks}
\def\checkendcite#1{\ifx\endcitedef#1?\else
	\expandafter\lookcite\expandafter#1\fi}
\def\lookcite#1 {\advance\citation by1\def\auxcite{#1}%
	\ifx\auxcite\curcite\the\citation\expandafter\endcitedef\else
	\expandafter\checkendcite\fi}
\def\cite#1{\makecite#1,\cite}
\def\makecite#1,#2{[\citenum{#1}\ifx\cite#2]\else\expandafter\clearcite\expandafter#2\fi}
\def\clearcite#1,\cite{, #1]}
\citedef
Be90
Bou
Cl83
Da
De
Doi84
Doi85
Doi92
Dor82
Ei
Ka
Kop77
Kop93
Kr81
Li03
Li
Ma91
Ma92
Ma94
Ma94b
MaD92
MaW94
Mc
Mo83
Mo
Mum
Ni89
Ni89b
Ni92
Par80
Rad77a
Rad77b
Row
Sch90
Sch92
Sch93
Sk93
Sw
Tak72
Tak77
Tak79
Zhu
\endcitedef
%
%
\def\references{\section References\par
	\bgroup
	\parindent=0pt
	\reffonts
	\rm
	\frenchspacing
	\setbox0\hbox{99. }\leftskip=\wd0
	}
\def\endreferences{\egroup}
\newtoks\nextauth
\newif\iffirstauth
\def\checkendauth#1{\ifx\endauth#1%
		\iffirstauth\the\nextauth
		\else{} and \the\nextauth\fi,
	\else\iffirstauth\the\nextauth\firstauthfalse
		\else, \the\nextauth\fi
		\expandafter\auth\expandafter#1\fi
	}
\def\auth#1,#2;{\nextauth={#1 #2}\checkendauth}
\newif\ifbookref
\def\nextref#1 {\par\hskip-\leftskip
	\hbox to\leftskip{\hfil\citenum{#1}.\ }%
	\initnextref}
\def\initnextref{\bookreffalse\firstauthtrue\ignorespaces}
\def\paper#1{{\it#1},}
\def\book#1{\bookreftrue{\it#1},}
\def\journal#1{#1}
\def\bookseries#1{#1,}
\def\Vol#1{\ifbookref Vol. #1,\else{\bf#1}\fi}
\def\publisher#1{#1,}
\def\Year#1{\ifbookref #1.\else(#1)\fi}
\def\Pages#1{#1.}
%
%
\newsymbol\smallsetminus 2272
\newsymbol\varnothing 203F
\newsymbol\square 1003
\newsymbol\rightsquigarrow 1320
\newsymbol\smallsetminus 2272
\newsymbol\varnothing 203F
\newsymbol\square 1003

\def\C#1{{\rm(C)}\ifcat a#1 \else\fi#1}
\def\eh{\hat e}
\def\F{{\cal F}}
\def\G{{\cal G}}
\def\g{\hbox{}\strut^g\!}
\def\I#1{{\cal I}_{#1}}
\def\mh{\hat m}
\def\M{{\cal M}}
\def\ADM{\hbox{}^{D'}_A\!\!\M}
\def\ADMH{\ADM^H}
\def\AGM{\hbox{}_{A,G}\M}
\def\AH{\HM_A}
\def\AHM{\hbox{}_{A\#H}\M}
\def\AM{\hbox{}_A\M}
\def\AMH{\AM^H}
\def\BM{\hbox{}_B\M}
\def\CM{\hbox{}^C\!\M}
\def\DM{\hbox{}^D\!\M}
\def\DMA{\DM_A}
\def\DAM{\hbox{}^D_A\!\M}
\def\DMAH{\hbox{}^D\!\MAH}
\def\DAMH{\hbox{}^D_A\!\M^H}
\def\HM{\hbox{}_H\M}
\def\MAH{\M_A^H}
\def\MD{\M^{D'}}
\def\SM{\hbox{}^S\!\M}
\def\SMA{\SM_A}
\def\SMAH{\hbox{}^S\!\MAH}

\def\P{{\cal C}_A(P)}
\def\U{U^{(X)}}
\def\V{{\cal V}}
\def\Z{{\Bbb Z}}
\def\l{\mathop{\rm lng}\nolimits}

\def\sq{\mathbin{\square}}
\def\sqH{\sq_{D'}H}
\def\limdir{\mathop{\vtop{\offinterlineskip\halign{##\hskip0pt\cr\rm lim\cr
	\noalign{\vskip1pt}
	$\scriptstyle\mathord-\mskip-10mu plus1fil
	\mathord-\mskip-10mu plus1fil
	\mathord\rightarrow$\cr}}}\nolimits}
\let\ot\otimes
\let\sbs\subset
\let\<\langle
\let\>\rangle
\def\Ann{\mathop{\rm Ann}\nolimits}
\def\can{{\rm can}}
\def\chr{\mathop{\rm char}\nolimits}
\def\Comm{Comm}
\def\cop{^{\mathop{\rm cop}}}
\def\End{\mathop{\rm End}\nolimits}
\def\ev{{\rm ev}}
\def\Fitt{\mathop{\rm Fitt}\nolimits}
\def\Hom{\mathop{\rm Hom}\nolimits}
\def\id{{\rm id}}
\def\Im{\mathop{\rm Im}}
\def\Ker{\mathop{\rm Ker}}
\def\mat#1{\left[\matrix{#1}\right]}
\def\Mat{\mathop{\rm Mat}\nolimits}
\def\Max{\mathop{\rm Max}\nolimits}
\def\op{^{\mathop{\rm op}}}
\def\rk{\mathop{\rm rk}\nolimits}
\def\soc{\mathop{\rm soc}\nolimits}
\def\Sp{\mathop{\frak S\frak p}\nolimits}
\def\Spec{\mathop{\rm Spec}}
\def\mapr#1{{}\mathrel{\smash{\mathop{\longrightarrow}\limits^{#1}}}{}}
\def\lmapr#1#2{{}\mathrel{\smash{\mathop{\count0=#1
  \loop
    \ifnum\count0>0
    \advance\count0 by-1\smash{\mathord-}\mkern-4mu
  \repeat
  \mathord\rightarrow}\limits^{#2}}}{}}
\def\mapd#1#2{\llap{$\vcenter{\hbox{$\scriptstyle{#1}$}}$}\big\downarrow
  \rlap{$\vcenter{\hbox{$\scriptstyle{#2}$}}$}}
\def\diagram#1{\vbox{\halign{&\hfil$##$\hfil\cr #1}}}
\def\diagramskip{\noalign{\smallskip}}
\let\al\alpha
\let\de\delta
\let\ep\varepsilon
\let\io\iota
\let\la\lambda
\let\ph\varphi
\let\th\theta
\let\ze\zeta
\let\Ga\Gamma
\let\De\Delta
\let\Om\Omega
\let\Th\Theta
%
%
\frontmatter

\title{Projectivity and freeness over comodule algebras}
\author{Serge Skryabin}
\address{Chebotarev Research Institute,
Universitetskaya St.~17, 420008
Kazan, Russia}
\email{Serge.Skryabin@ksu.ru}
\thanks{This research was supported by the project ``Construction and
applications of non-commutative geometry" from FWO Vlaanderen. I would like to
thank the Free University of Brussels VUB for hospitality during the time when
the work was conducted.}

\endfrontmatter

\abstract
Let $H$ be a Hopf algebra and $A$ an $H\!$-simple right $H\!$-comodule algebra.
It is shown that under certain hypotheses every $(H,A)$-Hopf module is either
projective or free as an $A$-module and $A$ is either a quasi-Frobenius or
a semisimple ring. As an application it is proved that every weakly finite
(in particular, every finite dimensional) Hopf algebra is free both as a left
and a right module over its finite dimensional right coideal subalgebras, and
the latter are Frobenius algebras. Similar results are obtained for
$H\!$-simple $H\!$-module algebras.
\endabstract

\section
Introduction

Let $H$ be a Hopf algebra over a field $k$. Starting from the work of Radford
\cite{Rad77a}, \cite{Rad77b} the question about the freeness and projectivity
of $H$ over its Hopf subalgebras aroused a substantial amount of interest.
Besides the pointed case studied in \cite{Rad77a} it is known that $H$ is a
projective module over any Hopf subalgebra whenever $H$ is commutative
\cite{Tak79, Cor.~1}. Takeuchi's original proof of this fact was based on the
faithful flatness as a preliminary step. By itself the faithful flatness of
$H$ over Hopf subalgebras is not so straightforward. When the subalgebra is
reduced one can apply the theorem on generic flatness; in general one has to
make a reduction to that case.

The first contribution to the theory that I am going to propose consists in
showing that the projectivity result just mentioned can be derived directly
from a known projectivity criterion in terms of the Fitting invariants. The
argument is very short and applies actually in a more general situation.
Suppose that $G$ is any covariant functor from the category of commutative
$k$-algebras to the category of groups. It makes sense to say what does it
mean for $G$ to operate on a commutative algebra $A$ by automorphisms
\cite{De}. If such an action is given one can introduce the notion of
$A,G$-modules. These are $A$-modules equipped with a compatible $G$-module
structure. It is shown in Proposition 1.1 that all Fitting invariants of an
$A$-finite $A,G$-module $M$ are $G$-stable ideals of $A$. If $A$ is $G$-simple
in the sense that $A$ has no nontrivial $G$-stable ideals, this immediately
implies the projectivity of $M$ as an $A$-module (Corollary 1.3 to Theorem
1.2). This gives a better result even in the case of an algebraic group $G$
operating rationally on $A$. Doraiswamy \cite{Dor82} treated connected
algebraic groups $G$ and finitely generated algebras $A$, but the method there
was based on the conclusion that $A$ is an integral domain under the
hypotheses stated. There are two cases where our result can be translated into
the language of Hopf algebras (Corollaries 1.5, 1.6). One can take $G$ to be
either a group scheme or a formal group scheme, and such group functors
correspond, respectively, to commutative and to cocommutative Hopf algebras.
In the first of these cases one obtains easily an improvement by removing the
finiteness condition on modules.

A generalization to noncommutative Hopf algebras presents a serious
difficulty. The celebrated achievement of Nichols and Zoeller consisted in
proving that every finite dimensional Hopf algebra $H$ is a free module over
its Hopf subalgebras \cite{Ni89}. In \cite{Ma92} Masuoka extended this result
by showing that $H$ is free both as a left and a right module over its right
coideal subalgebra $A$ if and only if $A$ is Frobenius. Several other
conditions equivalent to $A$ being Frobenius were given in \cite{Kop93},
\cite{Ma92}, \cite{MaD92}. The notion of coideal subalgebras appears to be of
fundamental importance. According to \cite{Ma92} there is a bijective
correspondence between the Frobenius right coideal subalgebras in $H$ and in
$H^*$. Whether all coideal subalgebras of $H$ are Frobenius was known to be
true under the assumption that $H$ has a cocommutative coradical \cite{Ma91}
and under the assumption that $H$ is involutory and $\chr k$ is either $0$ or
$>\!\dim H$ \cite{Ma92}. The primary motivation behind the present article
was to solve the last question for every finite dimensional $H$. Our approach
provides also a new proof of Nichols and Zoeller's theorem.

Again it is natural to work in more general settings. Let $H$ be an arbitrary
Hopf algebra over a field $k$ and $A$ a right $H\!$-comodule algebra. There are
categories of Hopf modules $\MAH$ and $\AMH$ introduced by Doi (e.g.,
\cite{Doi84}, \cite{Doi85}, \cite{Doi92}) as a generalization of their special
cases due to Takeuchi \cite{Tak72}, \cite{Tak79}. An $(H,A)$-Hopf module is
either right or left $A$-module equipped with a compatible $H\!$-comodule
structure. We say that $A$ is $H\!$-simple if $A$ has no nontrivial
$H\!$-costable (two-sided) ideals. Theorem 3.5 states that all objects of $\MAH$
are projective $A$-modules provided that $A$ is $H\!$-simple, semilocal and
satisfies one technical condition \C concerned with the weak finiteness of
certain rings. The last condition is not restrictive for many applications as,
for instance, left or right Noetherian rings, as well as rings finitely
generated as modules over commutative subrings are always weakly finite.
Condition \C is used to check that certain ideals of $A$ are $H\!$-costable. If
$\dim A<\infty$, in addition to previous hypotheses, then we will see in
Theorem 4.2 that $A$ is Frobenius and all objects of both $\MAH$ and $\AMH$
are projective $A$-modules. Theorem 4.5 ensures that $A$ is a quasi-Frobenius
ring under weaker assumptions about $A$. If $\dim H<\infty$, then the
hypotheses can be further weakened, and we also show in Theorem 5.2 that
$A$ is semisimple provided that so is $H$.

In order to apply the previous results to coideal subalgebras one has to know
that they are $H\!$-simple. This is based on Proposition 3.7: a right Artinian
$H\!$-comodule algebra $A$ satisfying \C is $H\!$-simple provided that there
exists a maximal ideal of $A$ containing no nonzero $H\!$-costable ideals.
Suppose next that $A$ is a finite dimensional right coideal subalgebra of $H$.
Under the assumption that $H$ is weakly finite the conclusions of Theorem 4.2
hold true for $A$. Thus $A$ is Frobenius not only for finite dimensional $H$
but also under much weaker assumptions. Moreover, all objects of both $\MAH$
and $\AMH$ are free $A$-modules, and the two categories are equivalent to the
categories of comodules over certain quotient coalgebras of $H$. In addition
the so-called normal basis property is fulfilled for $H$. All this is stated
in Theorem 6.1.

It is worth noting that the finiteness assumptions in the last result are
necessary. Nichols and Zoeller \cite{Ni89b} gave an example of a Hopf algebra
$H$ containing a two-dimensional Hopf subalgebra $A$ such that there exists a
finite dimensional object of $\AMH$ which is not a free $A$-module. We recall
this example in section 6 and show directly that $H$ is not weakly finite
here. On the other hand, the infinite dimensional coideal subalgebras are not
always $H\!$-simple and the projectivity over them may fail even when $H$ is
commutative. As an example take $H$ to be the group algebra of the free cyclic
group, say with a generator $g$, and $A$ the subalgebra of $H$ generated by
$g$. At least two related results were established in the literature without
finiteness restrictions on $H$, however. If $A$ is a finite dimensional Hopf
subalgebra of any $H$ and $A$ is either semisimple or normal in $H$, then $H$
is a free $A$-module \cite{Ni92}, \cite{Sch93}.

In section 7 of the paper we dualize the projectivity result to the
case of semilocal $H\!$-module algebras.
In conclusion I would like to thank A.~Masuoka and S.~Montgomery for helpful
comments.

\section
Notations and conventions

Let $k$ be the ground field. All algebras and coalgebras are over $k$, and $k$
serves normally as the base ring for functors $\ot$ and $\Hom$. If $A$ is an
algebra and $D$ a coalgebra, denote by $\AM$, $\M_A$, $\DM$, $\M^D$ the
categories of left $A$-modules, right $A$-modules, left $D$-comodules and
right $D$-comodules, respectively. Objects of the category $\DMA$ are vector
spaces equipped with a pair of commuting structures of a right $A$-module and
a left $D$-comodule so that all elements of $A$ operate as $D$-comodule
endomorphisms. The category $\DAM$ is defined similarly using left
$A$-modules.

Let $H$ be a Hopf algebra with the comultiplication $\De$, the counit $\ep$
and the antipode $s$. Either \cite{Mo} or \cite{Sw} can be used as a general
reference on Hopf algebras. A {\it right $H\!$-comodule algebra} is an algebra
$A$ together with a right $H\!$-comodule structure $\rho_A:A\to A\ot H$ such
that $\rho_A$ is a homomorphism of (unital) algebras. A {\it left $H\!$-module
algebra} is an algebra $A$ together with a left $H\!$-module structure such
that the map $\tau_A:A\to\Hom(H,A)$ defined by the rule $\tau_A(a)(h)=ha$ for
$a\in A$ and $h\in H$ is an algebra homomorphism. Here $\Hom(H,A)$ is
regarded as an algebra with respect to the convolution multiplication. We omit
the prefix ``right" for comodule algebras and the prefix ``left" for module
algebras. The subalgebra of invariants of an $H\!$-comodule algebra $A$ is
defined to be
$$
A^H=\{a\in A\mid\rho(a)=a\ot1\}.
$$

With each comodule algebra $A$ one associates the categories $\MAH$ and
$\AMH$. Their objects are either right or left $A$-modules together with a
compatible right $H\!$-comodule structure. For each module algebra $A$ we denote
by $\AH$ the category whose objects are right $A$-modules together with a
compatible left $H\!$-module structure. The compatibility condition in each of
the respective cases is as follows:
$$
\rho_M(ma)=\rho_M(m)\rho_A(a),\qquad
\rho_M(am)=\rho_A(a)\rho_M(m),
$$\removelastskip
$$
\tau_M(ma)=\tau_M(m)\tau_A(a)
$$
where $m\in M$ and $a\in A$. Here $\rho_M:M\to M\ot H$ is the comodule
structure map, and we regard $M\ot H$ as either right or left $A\ot H$-module
letting $A$ operate on the first tensorand and $H$ on the second tensorand via
multiplications. The map $\tau_M:M\to\Hom(H,M)$ is defined by the rule
$\tau_M(m)(h)=hm$ for $m\in M$ and $h\in H$, and we regard $\Hom(H,M)$ as a
right $\Hom(H,A)$-module with respect to the convolution action. Given linear
functions $\xi:H\to M$ and $\eta:H\to A$, one obtains $\xi\eta$ as the
composite
$$
H\mapr\De H\ot H\lmapr3{\xi\ot\eta}M\ot A\mapr{}M.
$$
We often omit the subscripts in the notations $\rho_A,\rho_M,\tau_A,\tau_M$.
Another way to express the compatibility of module and comodule structures in
the categories above is to say that the module structure map $M\ot A\to M$ or
$A\ot M\to M$ is a morphism in either $\M^H$ or $\HM$. Here we use the tensor
product of two comodule or two module structures defined in \cite{Mo, \S1.8}.
When necessary, we regard $k$ as a trivial $H\!$-module or $H\!$-comodule.

Objects of $\AH$ can be identified with the left modules over the smash
product algebra $A\op\#H\cop$. Here $A\op$ is $A$ taken with the opposite
multiplication and $H\cop$ is $H$ taken with the opposite comultiplication and
the same multiplication. If $A$ is commutative and $H$ cocommutative, then
$\AH\approx\AHM$.

An object $M$ of $\MAH$, $\AMH$, or $\AH$ will be called {\it$A$-finite} if it
is finitely generated as an $A$-module. We say that $M$ is {\it locally
$A$-finite} if it coincides with the union of its $A$-finite subobjects. The
same terminology will be used in the category $\AGM$ introduced in section 1.
Every object $M$ of either $\MAH$ or $\AMH$ is locally $A$-finite. Indeed,
the $A$-modules generated by the finite dimensional $H\!$-subcomodules of $M$
are $A$-finite subobjects of $M$.

An ``ideal" will be understood as a two-sided ideal unless explicitly
specified otherwise. The ideals in a (co)module algebra $A$ which are
respected by the (co)module structure will be termed $H\!$-(co)stable. The
algebra $A$ will be called {\it$H\!$-simple} if it has no nonzero proper
$H\!$-(co)stable ideals. Denote by $\Max R$ the set of all maximal ideals in a
ring $R$.

Note that $H$ is an $H\!$-comodule algebra with respect to $\De$ and $H$ is a
simple object of $\M_H^H$. By \cite{Sw, Th.~4.1.1} every $M\in\M_H^H$
decomposes as $M_0\ot H$ where $M_0=\{m\in M\mid\rho(m)=m\ot1\}$. For $M=H$
one has $M_0=k$ so that there is no room for nontrivial subobjects.

We use Sweedler's symbolic notations for comultiplication. If $h\in H$
and $m\in M$ where $M\in\M^H$, then
$$
\De(h)=\sum_{(h)}h_{(1)}\ot h_{(2)},\qquad
(\De\ot\id)\circ\De(h)=\sum_{(h)}h_{(1)}\ot h_{(2)}\ot h_{(3)},
$$\removelastskip
$$
\rho(m)=\sum_{(m)}m_{(0)}\ot m_{(1)},\qquad
(\rho\ot\id)\circ\rho(m)=\sum_{(m)}m_{(0)}\ot m_{(1)}\ot m_{(2)}.
$$

\section
1. Fully commutative case

Suppose that $A$ is a commutative ring, $M$ a finitely generated $A$-module.
For each $i\ge0$ the $i$th {\it Fitting invariant $\Fitt_i(M)$} of $M$ is
defined as follows. Taking an epimorphism of $A$-modules $\pi:F\to M$ where
$F$ is a free $A$-module, say of rank $n$, one sets $\Fitt_i(M)=A$ when $i\ge
n$; otherwise $\Fitt_i(M)$ is the ideal of $A$ generated by the determinants
of all $(n-i)\times(n-i)$ matrices $[f_j(x_l)]_{1\le j,l\le n-i}$ where
$f_1,\ldots,f_{n-i}$ run through $\Hom_A(F,A)$ and $x_1,\ldots,x_{n-i}$ run
through $\Ker\pi$. It is well known that this definition does not depend on
the choice of a presentation of $M$. We put formally $\Fitt_{-1}(M)=0$. We will
need two properties of the Fitting ideals (see \cite{Ei, Cor. 20.5 and Prop.
20.8}):

\jtem(F1)
$\Fitt_i(B\ot_AM)=\Fitt_i(M)B$ whenever $B$ is a commutative $A$-algebra,

\jtem(F2)
For $M$ to be a projective $A$-module of constant rank $r\ge0$ it is
necessary and sufficient that $\Fitt_r(M)=A$ and $\Fitt_{r-1}(M)=0$.

Denote by $\Comm_k$ the category of commutative $k$-algebras. A {\it group
$k$-functor} $G$ is any functor from $\Comm_k$ to the category of groups
\cite{De}. Thus $G$ associates a group $G(R)$ with each $R\in\Comm_k$ and a
group homomorphism $G(\ph):G(R)\to G(R')$ with each homomorphism of
commutative algebras $\ph:R\to R'$. A {\it$G$-module\/} is a vector space
$V$ together with $R$-linear actions of the groups $G(R)$ on the $R$-modules
$V\ot R$ which are given for each $R\in\Comm_k$ and are compatible with
morphisms in $\Comm_k$, that is, whenever $g\in G(R)$ and $\ph:R\to R'$ is an
algebra homomorphism, the transformation of $V\ot R'$ afforded by
$G(\ph)(g)\in G(R')$ is the $R'$-linear extension of the transformation of
$V\ot R$ afforded by $g$ \cite{De, Ch.~II, \S2}.

Assume further that $A\in\Comm_k$. One says that a group $k$-functor $G$ {\it
operates on $A$ by automorphisms} if $A$ is given a $G$-module structure and
for each $R\in\Comm_k$ the group $G(R)$ acts on $A\ot R$ as a group of algebra
automorphisms. If such an action is given, an ideal $I$ of $A$ will be called
{\it$G$-stable} if $I\ot R$ is stable under $G(R)$ for each $R\in\Comm_k$. The
definition of $G$-stable subalgebras is similar. We say that $A$ is
{\it$G$-simple} if $A$ has no nonzero proper $G$-stable ideals.

An {\it$A,G$-module} $M$ is an $A$-module together with a $G$-module structure
which satisfies the following compatibility condition:
$$
g(am)=(ga)(gm)\qquad{\rm for\ all\ }g\in G(R),\ a\in A\ot R,\ m\in M\ot R.
\eqno(*)
$$
Denote by $\AGM$ the category of all $A,G$-modules.

If $G$ is a group scheme, then the notion of $A,G$-modules can be interpreted
geometrically in terms of $G$-linearized quasicoherent sheaves on the affine
scheme $\Spec A$ \cite{Mum, Ch.~I, \S3}. Such structures have been studied
also in purely algebraic context. For instance, if $\Ga$ is an ordinary group
and $G$ the constant group functor such that $G(R)=\Ga$ for each $R\in\Comm_k$
and $G(\ph)$ is the identity map $\Ga\to\Ga$ for each morphism $\ph$ in
$\Comm_k$, an $A,G$-module is just a module over the skew group ring $A*\Ga$.

\proclaim
Proposition 1.1.
If $M\in\AGM$ is $A$-finite, then all its Fitting ideals $\,\Fitt_i(M)$ are
$G$-stable.
\endproclaim

\Proof.
Suppose that $g\in G(R)$ where $R\in\Comm_k$. For every $A\ot R$-module $N$
one can define a new $A\ot R$-module $\g N$ such that $\g N=N$ as abelian
groups and each $a\in A\ot R$ operates in $\g N$ as $g^{-1}(a)$ does in $N$.
We claim that
$$
\Fitt_i(\g N)=g\bigl(\Fitt_i(N)\bigr).\eqno(**)
$$
To prove this equality consider an epimorphism of $A\ot R$-modules $\pi:F\to
N$ where $F$ is a free $A\ot R$-module, say of rank $n$. The same map $\pi$ is
also an epimorphism of $A\ot R$-modules $\g F\to\g N$. Clearly $\g F$ is a
free $A\ot R$-module of rank $n$. Indeed, any basis for $F$ is also a basis
for $\g F$. In particular, both sides of $(**)$ are equal to $A\ot R$ when
$i\ge n$. Assume that $0\le i<n$. There is a bijection
$$
\Hom_{A\ot R}(F,A\ot R)\to\Hom_{A\ot R}(\g F,A\ot R)
$$
given by $f\mapsto\g f$ where $(\g f)(x)=g(f(x))$ for $x\in F$. As $g$ operates
on $A\ot R$ as an algebra automorphism, we have
$$
\det[(\g f_j)(x_l)]_{1\le j,l\le n-i}=g\bigl(\det[f_j(x_l)]_{1\le j,l\le n-i}\bigr)
$$
whenever $f_1,\ldots,f_{n-i}\in\Hom_{A\ot R}(F,A\ot R)$ and
$x_1,\ldots,x_{n-i}\in\Ker\pi$. The determinants on the left and right hand
sides of this formula generate the two ideals in $(**)$, whence the claim.

We will apply $(**)$ to the $A\ot R$-module $N=M\ot R$. First of all, $\g
N\cong N$ in this case because the transformation of $N$ afforded by $g$ is an
isomorphism between the two $A\ot R$-module structures according to $(*)$.
Hence $\Fitt_i(\g N)=\Fitt_i(N)$ so that $(**)$ just says that $\Fitt_i(N)$ is
stable under $g$. Property (F1) applied to the $A$-algebra $B=A\ot R$ shows,
however, that $\Fitt_i(N)=\Fitt_i(M)\ot R$.
\endproof

\proclaim
Theorem 1.2.
Suppose that $B$ is a commutative $A$-algebra such that $IB=B$ for every
$G$-stable ideal $I\ne0$ of $A$. If $M\in\AGM$ is locally $A$-finite, then
$B\ot_AM$ is a projective $B$-module. The latter has constant rank
whenever $M$ is $A$-finite.
\endproclaim

\Proof.
Suppose first that $M$ is $A$-finite. Since $\Fitt_i(M)$ is a $G$-stable ideal
of $A$ by Proposition 1.1, it follows from (F1) that $\Fitt_i(B\ot_AM)$ is
either $0$ or $B$ for each $i\ge-1$. On the other hand, the latter ideal is
$0$ for $i=-1$ and is $B$ for sufficiently large $i$. There exists then
$r\ge0$ such that $\Fitt_i(B\ot_AM)$ is $B$ for $i=r$ and is $0$ for $i=r-1$.
By (F2) the $B$-module $B\ot_AM$ is projective of rank $r$.

Consider now the general case. Denote by $\G$ and $\F$ the sets whose elements
are all $A,G$-submodules of $M$ and the $A$-finite ones, respectively. Given
$N,N'\in\G$ such that $N\sbs N'$, put
$$
T_{NN'}=\Im(B\ot_AN\lmapr3\can B\ot_AN').
$$

{\it Step\/ }1.
Put also $T_N=T_{NM}$ for short. Thus $T_N\sbs B\ot_AM$ is a $B$-submodule.
We will show that $T_N$ is a projective $B$-module whenever $N\in\F$. The
set $\F$ is directed by inclusion and $M=\bigcup_{N'\in\F}N'$ by local
finiteness of $M$. Since tensor products commute with direct limits, we
have $B\ot_AM\cong\limdir_{N'\in\F}B\ot_AN'$ and
$$
T_N\cong\limdir_{N'\in\F_N}T_{NN'}
$$
where $\F_N=\{N'\in\F\mid N\sbs N'\}$. If $N'\in\F_N$ then both $N'$ and
$N'/N$ are $A$-finite objects of $\AGM$, whence $B\ot_AN'$ and $B\ot_AN'/N$
are projective $B$-modules of constant rank. Since $(B\ot_AN')/T_{NN'}\cong
B\ot_AN'/N$ by right exactness of tensor products, we deduce that $T_{NN'}$ is
a projective $B$-module of constant rank. Suppose that $N',N''\in\F_N$ and
$N'\sbs N''$. Then the canonical morphism $B\ot_AN'\to B\ot_AN''$ in $\BM$
maps $T_{NN'}$ onto $T_{NN''}$. The kernel, say $K$, of the induced
epimorphism $T_{NN'}\to T_{NN''}$ is a projective $B$-module of constant rank
equal to $\rk T_{NN'}-\rk T_{NN''}$. We see that $\rk T_{NN'}\ge\rk T_{NN''}$
and one has an equality here if and only if $K=0$. Pick $N'\in\F_N$ for which
$\rk T_{NN'}$ attains the minimum value. Then $T_{NN'}$ is mapped
isomorphically onto $T_{NN''}$ for each $N''\in\F_N$ such that $N'\sbs N''$,
and it follows that $T_N\cong T_{NN'}$.

{\it Step\/ }2.
Next we will prove that $T_{N'}/T_N$ is a projective $B$-module whenever
$N,N'\in\G$ are such that $N\sbs N'$ and $N'/N$ is $A$-finite. Since tensor
products are right exact, we have $(B\ot_AM)/T_N\cong B\ot_AM/N$. This induces
an isomorphism of $T_{N'}/T_N$ onto the image of the canonical map
$B\ot_AN'/N\to B\ot_AM/N$. It remains to apply Step 1 to the locally
$A$-finite object $M/N\in\AGM$.

{\it Final Step}.
Let $\eta:K\to L$ be an epimorphism and $\xi:B\ot_AM\to L$ any morphism in
$\BM$. We have to find a morphism $\ze:B\ot_AM\to K$ in $\BM$ such that
$\xi=\eta\circ\ze$. Consider the set $\Om$ of all pairs $(N,\ze)$ where
$N\in\G$ and $\ze:T_N\to K$ is a morphism in $\BM$ such that
$\eta\circ\ze=\xi|_{T_N}$. Note that $(0,0)\in\Om$. For two pairs in $\Om$ set
$(N,\ze)\le(N',\ze')$ if and only if $N\sbs N'$ and $\ze=\ze'|_{T_N}$. By
Zorn's lemma $\Om$ has a maximal element. Let now $(N,\ze)$ be such a maximal
element. Suppose $N\ne M$. Then there exists $F\in\F$ such that $F\not\sbs N$.
Put $N'=N+F$. By Step 2 $T_{N'}=V\oplus T_N$ for some projective $B$-submodule
$V$. The restriction of $\xi$ to $V$ factors through $K$ in $\BM$, and this
shows that $(N',\ze')\in\Om$ for a suitable extension $\ze'$ of $\ze$. We have
obtained a contradiction with the maximality of $(N,\ze)$. Thus $N=M$, whence
the required $\ze$.
\endproof

\proclaim
Corollary 1.3.
Let $M\in\AGM$ be locally $A$-finite. If $P$ is a prime ideal of $A$
containing no nonzero $G$-stable ideals, then the localization $M_P$ of $M$ at
$P$ is a free module over the local ring $A_P$. If $A$ is $G$-simple, then
$M$ is a projective $A$-module.
\endproclaim

\Proof.
Apply Theorem 1.2 taking $B=A_P$ in one case and $B=A$ in the other. By a well
known theorem of Kaplansky all projective modules over a local ring are free
\cite{Bou, Ch.~II, \S3, Exercise 3}.
\endproof

\proclaim
Corollary 1.4.
Suppose $B$ is a $G$-simple commutative algebra, $A\sbs B$ a $G$-stable
subalgebra. If $B$ is locally $A$-finite as an object of $\AGM$, then $B$ is
flat over $A$.
\endproclaim

\Proof.
If $I$ is any nonzero $G$-stable ideal of $A$ then $IB$ is a nonzero
$G$-stable ideal of $B$, whence $IB=B$. Suppose that $Q$ is a prime ideal
of $B$ and $P=Q\bigcap A$. Then $P$ can contain no nonzero $G$-stable ideals
of $A$, and so $B_P$ is a free $A_P$-module by Corollary 1.3. Then $B_Q$ is
flat over $A_P$. The flatness of $B$ follows now from \cite{Bou, Ch.~II, \S3,
Prop.~15}.
\endproof

{\it Examples.}
1) Every cocommutative Hopf algebra $H$ determines a formal group scheme
$G=\Sp^*\!H$ which can be regarded as a group $k$-functor such that $G(R)$ is
the group of all grouplike elements of the Hopf algebra $H\ot R$ over $R$.
Thus
$$
G(R)=\{g\in H\ot R\mid\De_R(g)=g\ot g\ {\rm and}\ \ep_R(g)=1\}
\qquad{\rm for\ }R\in\Comm_k
$$
where $\De_R$ and $\ep_R$ are the comultiplication and the counit of $H\ot R$.
We follow \cite{Tak77} in considering a formal scheme as a functor defined on
the whole $\Comm_k$ rather than on the full subcategory consisting of finite
dimensional commutative algebras. In this case there is a bijective
correspondence between the $G$-module and $H\!$-module structures. To let $G$
operate on $A$ by automorphisms is the same as to make $A$ into an $H\!$-module
algebra, and one has $\AGM\approx\AHM$. A special case of Theorem 1.2 thus
yields

\proclaim
Corollary 1.5.
Suppose that $H$ is a cocommutative Hopf algebra, $A$ a commutative $H\!$-module
algebra, and $B$ a commutative $A$-algebra such that $IB=B$ for every
$H\!$-stable ideal $I\ne0$ of $A$. If $M\in\AHM$ is locally $A$-finite, then
$B\ot_AM$ is a projective $B$-module.
\endproclaim

In particular, one may take $H$ to be a group algebra. In another case where
$H$ is the universal enveloping algebra of a Lie algebra one obtains the
projectivity result stated in \cite{Sk93, Th.~1.6}.

\smallskip
2) Every commutative Hopf algebra $H$ determines an affine group scheme
$G=\Sp H$ which can be regarded as a group $k$-functor such that $G(R)$ is
the set of all algebra homomorphisms $H\to R$. The multiplication in $G(R)$
comes from the convolution multiplication in $\Hom(H,R)$. The $G$-modules here
are precisely right $H\!$-comodules \cite{De, Ch.~II, \S2, n$^\circ$2.1}. Giving
an action of $G$ on $A$ by automorphisms makes $A$ into an $H\!$-comodule
algebra and vice versa. There is a category equivalence $\AGM\approx\MAH$.

\proclaim
Corollary 1.6.
Suppose that $H$ is a commutative Hopf algebra, $A$ a commutative $H\!$-comodule
algebra, and $B$ a commutative $A$-algebra such that $IB=B$ for every
$H\!$-costable ideal $I\ne0$ of $A$. Then $B\ot_AM$ is a projective $B$-module
for any $M\in\M_A^H$.
\endproclaim

We know that $H$ is an $H\!$-simple $H\!$-comodule algebra, that is, a $G$-simple
algebra in this case. A subalgebra of $H$ is $G$-stable if and only if it is
$H\!$-costable. Corollary 1.4 includes therefore the following result of
Masuoka and Wigner \cite{MaW94, Th.~3.4}: {\it a commutative Hopf algebra is a
flat module over every right coideal subalgebra}. This statement translates
into geometric language: {\it if an affine group scheme $G$ operates on an
affine scheme $X$, then every $G$-equivariant morphism $f:G\to X$ decomposes
as $G\to Y\to X$ where the first morphism is flat and the second one is a
closed immersion}. Indeed, if $X=\Sp A$, then $f$ corresponds to a
homomorphism of $H\!$-comodule algebras $\ph:A\to H$. The image of $\ph$ is a
right coideal subalgebra of $H$ over which $H$ is flat. We may take therefore
$Y=\Sp\ph(A)$. In a special case where $f$ is a homomorphism of affine group
schemes, $f$ induces a faithfully flat homomorphism of $G$ onto a closed group
subscheme of $X$ \cite{De, Ch.~III, \S3, Cor.~7.3}. A purely algebraic
formulation of the latter fact says that {\it a commutative Hopf algebra is a
faithfully flat module over every Hopf subalgebra} \cite{Tak72, Th.~3.1}. The
morphism $G\to Y$ in the decomposition of $f$ above is always faithfully flat
also under the assumption that $X$ is a finite scheme. Indeed, as was
announced by Masuoka \cite{Ma94b, Th.~3.5} (and follows also from Theorem 6.1
later in the paper), {\it a commutative Hopf algebra is a free module
over every finite dimensional right coideal subalgebra}. In case of Hopf
subalgebras the last statement was proved by Radford \cite{Rad77b}.

\section
2. Some ring-theoretic facts

Recall that a ring is {\it semilocal} if its factor ring by the Jacobson
radical is semisimple Artinian. Following \cite{Row, Def.~1.3.30} we say that
a ring $R$ is {\it weakly $n$-finite} if for every pair of $n\times n$-matrices
$X,Y\in\Mat_n(R)$ the equality $XY=1$ in the matrix ring implies the equality
$YX=1$. A ring $R$ is {\it weakly finite} (or {\it stably finite} according to 
different sources, e.g. \cite{Mo83}) if $R$ is weakly $n$-finite for all 
integers $n>0$.

\proclaim
Lemma 2.1.
Let $A$ be a ring, $M$ a finitely generated right $A$-module. The ring
$R=\End_AM$ is weakly finite in each of the following three cases{\rm:}

\item(a)
every finite subset of $A$ is contained in a right Noetherian subring{\rm,}

\item(b)
$A$ is commutative{\rm,}

\item(c)
$A$ is weakly finite and $M$ is projective.

\endproclaim

\Proof.
Let $\xi,\eta\in R$ be such that $\xi\circ\eta=\id$. Pick any generators
$e_1,\ldots,e_m$ for $M$ over $A$. Under hypothesis (a) there exists a
Noetherian subring $B\sbs A$ such that $\xi(e_1),\ldots,\xi(e_m)$ and
$\eta(e_1),\ldots,\eta(e_m)$ lie in the $B$-submodule $N\sbs M$ generated by
$e_1,\ldots,e_m$. Then $N$ is stable under both $\xi$ and $\eta$. In
particular, $\xi$ induces a surjective endomorphism of $N$. Since $N$ is a
Noetherian $B$-module, $\xi$ has to be bijective on $N$ \cite{Ka, Th.~6.4.1},
and it follows that $\eta\circ\xi$ is identity on $N$. Since $M=NA$,
we conclude that $\eta\circ\xi=\id$ on the whole $M$. Note that hypothesis (b)
is a special case of (a). Indeed, finitely generated commutative rings are all
Noetherian. Under hypothesis (c) $M$ is a direct summand of a free right
$A$-module $A^t$ for some integer $t\ge0$. Letting both $\xi$ and $\eta$ act
as identity on a complementary summand, we extend $\xi$ and $\eta$ to
endomorphisms of $A^t$ so that $\xi\circ\eta=\id$ on $A^t$. Since the ring
$\End_AA^t\cong\Mat_t(A)$ is weakly 1-finite, the equality $\eta\circ\xi=\id$
holds again.

We have checked that $R$ is weakly 1-finite in each of cases (a), (b), (c). It
remains to observe that $\Mat_n(R)\cong\End_AM^n$ and $M^n$ is a finitely
generated right $A$-module which is also projective in case (c).
\endproof

\proclaim
Proposition 2.2.
A ring $R$ is weakly finite in each of the cases listed below{\rm:}

\item(a)
$R$ is left or right Noetherian{\rm,}

\item(b)
$R$ is semilocal{\rm,}

\item(c)
$R\op$ is weakly finite{\rm,}

\item(d)
$R$ is a finitely generated left or right module over a commutative
subring $S${\rm,}

\item(e)
$R$ is a finitely generated projective left or right module over a weakly
finite subring $S${\rm,}

\item(f)
$R=A\ot B$ where $A$, $B$ are algebras over the ground field $k$ such that $A$
is a finitely generated module over its center $Z$ and $B$ is weakly finite.

\endproclaim

Part (f) is a special case of Montgomery's result \cite{Mo83, Th.~1} where the 
assumption on $A$ is weakened to $A$ being any polynomial identity algebra.
We nevertheless provide a proof of (f) as it uses a different argument.

\Proof.
Part (a) is covered by \cite{Mo83, Cor.~1} or \cite{Row, Th.~3.2.37}.

(b) Let $J$ be the Jacobson radical of $R$. The ring $R/J$ is Artinian, hence
weakly finite. Then $R$ is weakly finite too \cite{Mo83, Lemma 2}.

(c) Use an isomorphism $\Mat_n(R\op)\cong\Mat_n(R)\op$.

(d), (e) There are embeddings $R\sbs\End R_S$ and $R\sbs(\End\hbox{}_SR)\op$.
Now we can apply Lemma 2.1 and (c).

(f) Since $\Mat_n(R)\cong A\ot\Mat_n(B)$ and the ring $\Mat_n(B)$ is weakly
finite in view of (e), it suffices to prove that $R$ is weakly 1-finite. Let
$e_1,\ldots,e_m$ generate $A$ as a $Z$-module. There exists a finitely
generated subalgebra $Z'\sbs Z$ such that $e_ie_j\in\sum_{l=1}^mZ'e_l$ for all
$1\le i,j\le m$. Then $\sum_{l=1}^mKe_l$ is a subalgebra of $A$ whenever $K$ is
a subalgebra satisfying $Z'\sbs K\sbs Z$. This shows that every finite subset of
$A$ is contained in a subalgebra $A'\sbs A$ which is finitely generated as a
module over some finitely generated central subalgebra $K\sbs A'$. Given $x,y\in
R$, there exists $A'$ as above such that $x,y\in A'\ot B$. Suppose that
$xy=1$. If $I$ is any ideal of $K$ such that $\dim K/I<\infty$, then $\dim
A'/IA'<\infty$ as well. The ring $A'/IA'\ot B$ is weakly finite by (e) as it
is a finitely generated free module over its subring $1\ot B$. Considering
the images of $x,y$ in $A'/IA'\ot B$, we deduce that $yx-1\in IA'\ot B$. Hence
$yx-1$ lies in the intersection $\bigcap\,IA'\ot B$ over all ideals $I$ of $K$
satisfying $\dim K/I<\infty$. It remains to check that $\bigcap\,IA'=0$. If
$P\in\Max K$, then $\dim K/P<\infty$ by Hilbert's Nullstellensatz. As $K$ is
Noetherian, we have also $\dim K/P^i<\infty$ for all $i>0$. For an element
$a\in A'$ denote $\Ann(a)=\{z\in K\mid za=0\}$. By \cite{Bou, Ch.~III, \S3,
Prop.~5} the inclusion $a\in\bigcap_{i>0}P^iA'$ implies that $\Ann(a)\not\sbs
P$, and if this is valid for all $P\in\Max K$ then $\Ann(a)=K$, which is only
possible for $a=0$.
\endproof

Let $M$ be a finitely generated right $R$-module. For every system of its
generators $e_1,\ldots,e_n$ denote by $\I{e_1,\ldots,e_n}$ the ideal of $R$
generated by all elements of $R$ which occur as a coefficient in a zero linear
combination $e_1x_1+\cdots+e_nx_n=0$ with $x_1,\ldots,x_n\in R$. If $R$ is
commutative then clearly $\I{e_1,\ldots,e_n}=\Fitt_{n-1}(M)$.

\proclaim
Lemma 2.3.
Let $e_1,\ldots,e_n$ and $e'_1,\ldots,e'_n$ be two systems of generators for
$M$ having the same number of elements. Suppose that $I$ is an ideal of $R$
such that $\I{e_1,\ldots,e_n}\sbs I$. If the ring $R/I$ is weakly $n$-finite
then $\I{e'_1,\ldots,e'_n}\sbs I$.
\endproclaim

\Proof.
There are expressions $e'_j=\sum_{i=1}^ne_ia_{ij}$ and
$e_j=\sum_{i=1}^ne'_ib_{ij}$ for some $a_{ij},b_{ij}\in R$ with $1\le i,j\le
n$. Then $e_l=\sum_{j=1}^n\sum_{i=1}^ne_ia_{ij}b_{jl}$, i.e.
$$
\sum_{i=1}^n\,e_i\Bigl(\,\sum_{j=1}^na_{ij}b_{jl}-\de_{il}\Bigr)=0
$$
for each $l=1,\ldots,n$. This shows that $\sum_{j=1}^na_{ij}b_{jl}-\de_{il}\in
I$ for all $i,l$. Denote by $X,Y\in\Mat_n(R/I)$ the $n\times n$ matrices whose
entries are the cosets, respectively, of elements $a_{ij}$ and $b_{ij}$ modulo
$I$. Previous inclusions can be rewritten as a single matrix equality telling
us that $XY$ is the identity matrix in $\Mat_n(R/I)$. By the hypotheses of the
lemma so is then the product $YX$ as well, i.e.
$$
\sum_{j=1}^nb_{ij}a_{jl}-\de_{il}\in I
\qquad{\rm for\ all\ }1\le i,l\le n.\eqno(*)
$$
Suppose now that $\sum_{l=1}^ne'_ly_l=0$ for some $y_1,\ldots,y_n\in R$.
Expressing each $e'_l$ as a linear combination of $e_1,\ldots,e_n$, we get
$\sum_{j=1}^ne_j\bigl(\sum_{l=1}^na_{jl}y_l\bigr)=0$ so that
$\sum_{l=1}^na_{jl}y_l\in I$ for all $j=1,\ldots,n$. Since $I$ is a two-sided
ideal, preceding inclusions together with $(*)$ yield
$$
y_i=\sum_{j=1}^n\,b_{ij}\sum_{l=1}^na_{jl}y_l
-\sum_{l=1}^n\Bigl(\sum_{j=1}^nb_{ij}a_{jl}-\de_{il}\Bigr)y_l\in I
$$
for all $i=1,\ldots,n$. The inclusion asserted by the lemma follows now from
the definition of $\I{e'_1,\ldots,e'_n}$.
\endproof

\proclaim
Lemma 2.4.
Let $R$ be a semilocal ring and $\V$ the class of all finitely generated
projective right $R$-modules $V$ such that a direct sum of finitely many
copies of $V$ is a free $R$-module. In order that $V\in\V$ be a free
$R$-module it is necessary and sufficient that $V/VQ$ be a free $R/Q$-module
for at least one $Q\in\Max R$. There exists $U\in\V$ such that all modules in
$\V$ are direct sums of copies of $U$.
\endproclaim

\Proof.
Denote by $J$ the Jacobson radical of $R$. Suppose that $V$ is a finitely
generated projective right $R$-module. We have
$V/V\!J\cong\prod_{Q\in\Max R}V/VQ$ and each $V/VQ$ is a direct sum of
finitely many copies of the single simple $R/Q$-module. Put
$$
r_Q(V)={\l V/VQ\over\l R/Q}\qquad{\rm for\ }Q\in\Max R
$$
where $\l$ is a shorthand for length. In order that $V$ be a free $R$-module
of rank $n$ it is necessary and sufficient that $V/V\!J$ be a free $R/J$-module
of rank $n$ \cite{Bou, Ch.~II, \S3, Prop.~5}. The latter condition means
precisely that each $V/VQ$ is a free $R/Q$-module of rank $n$ or,
equivalently, $r_Q(V)=n$ for each $Q$. If $V^t$ is a free $R$-module of rank
$m$ for some integers $t>0$ and $m\ge0$, then $t\l V/VQ=m\l R/Q$, and so
$r_Q(V)=m/t$ for all $Q\in\Max R$. This shows that the numbers $r_Q(V)$ for
$V\in\V$ do not depend on $Q$. Let $r(V)$ denote the common value of these
numbers. We see that $V\in\V$ is free if and only if $r(V)\in\Z$. This
proves the first assertion.

We claim that, whenever $r(V)\ge r(W)$ for some $V,W\in\V$, there always
exists an epimorphism $V\to W$. First of all, there exists an epimorphism
$V/VQ\to W/WQ$ for each $Q\in\Max R$, hence an epimorphism $V/V\!J\to W/W\!J$.
The latter can be lifted to a homomorphism $\xi:V\to W$ since $V$ is
projective. By Nakayama's lemma $\xi$ is surjective.

We can find an integer $d>0$ such that $r(V)\in{1\over d}\Z$ for all
$V\in\V$. For instance, we may take $d=\l R/Q$ for any chosen $Q\in\Max R$.
It follows that there exists $0\ne U\in\V$ such that $r(U)\le r(V)$ for all
$0\ne V\in\V$. Given $V$, let $n\ge0$ be the largest integer such that
$r(U^n\!)=r(U)n\le r(V)$. As we have seen, $U^n$ is a homomorphic image of $V$,
so that $V\cong U^n\oplus T$ for some projective $R$-module $T$. Now
$r_Q(T)=r(V)-r(U^n\!)$ for all $Q\in\Max R$, and so $r_Q(T^d\!)$ is an integer
not depending on $Q$. This shows that $T^d$ is a free $R$-module, whence
$T\in\V$. By the choice of $n$ we must have $r(T)<r(U)$, but then $T=0$ and
$V\cong U^n$.
\endproof

\proclaim
Lemma 2.5.
Let $R$ be a semilocal ring. A right $R$-module $M$ is necessarily free as long
as $M$ is not finitely generated and there exists a family $\F$ of its
submodules satisfying conditions {\rm(a)} and {\rm(b)} below{\rm:}

\item(a)
$0\in\F$ and the union of every chain in $\F$ is again in $\F${\rm,}

\item(b)
each $N\in\F\!$, $N\ne M$, is properly contained in some $N'\in\F$ with
$N'\!/N\in\V$.

\endproclaim

\Proof.
Let $U$ be as in Lemma 2.4. For each subset $X\sbs\Hom_R(U,M)$ denote by $\U$
the direct sum of the family of copies of $U$ indexed by $X$. For $\xi\in X$
let $\io_\xi:U\to\U$ be the canonical embedding of a summand. There exists a
unique homomorphism $\th_X:\U\to M$ such that $\th_X\circ\io_\xi=\xi$ for each
$\xi$.

Denote by $\Om$ the set of all pairs $(N,X)$ where $N\in\F$ and
$X\sbs\Hom_R(U,M)$ is a subset such that $\th_X$ is an isomorphism onto $N$.
Note that $(0,\varnothing)\in\Om$. Define a partial order on $\Om$ by setting
$(N,X)\le(N',X')$ for two pairs in $\Om$ if and only if $N\sbs N'$ and $X\sbs
X'$. If $\{(N_\al,X_\al)\}$ is a chain in $\Om$ indexed by elements $\al$ of
some set then $({}\bigcup N_\al,\bigcup X_\al)$ is in $\Om$ according to
(a). By Zorn's lemma $\Om$ has a maximal element. Let now $(N,X)$ be such a
maximal element. Suppose that $N\ne M$, and let $N'$ be as in (b). We have
then an isomorphism of $R$-modules $N'\cong N\oplus V$ where $V=N'\!/N\in\V$.
By Lemma 2.4 $V\cong U^n$ for some integer $n\ge0$. Let $\eta_1,\ldots,\eta_n$
be the isomorphisms of $U$ onto direct summands in the decomposition of $V$.
Denote $X'=X\bigcup{}\{\eta_1,\ldots,\eta_n\}$. Then $(N',X')\in\Om$ and
$(N,X)<(N',X')$. We have arrived at a contradiction with the maximality of
$(N,X)$.

It follows that any maximal element of $\Om$ is necessarily $(M,X)$
for some $X$. Thus $M\cong\U$. Since $M$ is not finitely generated, $X$ must
be infinite.  Then $X$ has the same cardinality as the set
$X\times\{1,\ldots,t\}$ for any integer $t>0$. However, $U^t$ is a free
$R$-module for a suitable $t$. It follows that $\U\cong R^{(X)}$.
\endproof

\section
3. Projectivity result for comodule algebras

Let $A$ be an $H\!$-comodule algebra. In the main result of this section we
encounter the following technical condition on $A$:

\jtem(C)
the ring $A/Q\ot H$ is weakly finite for each $Q\in\Max A$.

Using Proposition 2.2 one can verify the validity of this condition in many
cases. For instance, (C) always holds whenever $H$ is weakly finite and $A$ is
finitely generated as a module over its center. If $\dim H<\infty$, condition
(C) is fulfilled provided that $A$ is semilocal. The same is true if $H$ is
only {\it residually finite dimensional} in the sense that its ideals of
finite codimension have zero intersection.

\proclaim
Lemma 3.1.
If $e_1,\ldots,e_n$ generate $M\in\MAH$ as an $A$-module, then
$\rho(e_1),\ldots,\rho(e_n)$ generate $M\ot H$ as an $A\ot H$-module.
\endproclaim

\Proof.
If $m\in M$, then
$$
m\ot1=\sum_{(m)}\,m_{(0)}\ot m_{(1)}s(m_{(2)}\!)
=\sum_{(m)}\,\rho(m_{(0)}\!)\cdot\bigl(1\ot s(m_{(1)}\!)\bigr),
$$
which shows that $M\ot H=\rho(M)\cdot(A\ot H)$. Since
$\rho(ma)=\rho(m)\rho(a)$ for $m\in M$ and $a\in A$, we see that $\rho(M)$ is
generated by $\rho(e_1),\ldots,\rho(e_n)$ as a $\rho(A)$-module, and the lemma
is proved.
\endproof

\proclaim
Lemma 3.2.
Let $M\in\MAH$ be generated by $e_1,\ldots,e_n$ as an $A$-module. Suppose that
$I$ is an ideal of $A$ such that $\I{e_1,\ldots,e_n}\sbs I$ and the ring
$A/I\ot H$ is weakly $n$-finite. Then there exists an $H\!$-costable ideal
$K$ of $A$ such that $\I{e_1,\ldots,e_n}\sbs K\sbs I$.
\endproclaim

\Proof.
Put $J=\I{e_1,\ldots,e_n}$  (as defined in Lemma 2.3) and
$K=\rho^{-1}(I\ot H)$. Since $\rho:A\to A\ot H$ is an algebra homomorphism,
$K$ is an ideal of $A$. One has $K\sbs I$ because $(\id\ot\ep)\circ\rho=\id$
and $\rho(K)\sbs K\ot H$ because
$$
(\rho\ot\id)\circ\rho(K)=(\id\ot\De)\circ\rho(K)\sbs(\id\ot\De)(I\ot H)\sbs
I\ot H\ot H.
$$
The $A\ot H$-module $M\ot H$ is clearly generated by $e_1\ot1,\ldots,e_n\ot1$.
It is straightforward to see that $\I{e_1\ot1,\ldots,e_n\ot1}=J\ot H$.
By Lemma 3.1 $\rho(e_1),\ldots,\rho(e_n)$ is another system of generators for
the $A\ot H$-module $M\ot H$. Applying Lemma 2.3 with $R=A\ot H$ and the ideal
$I\ot H$, we get $\I{\rho(e_1),\ldots,\rho(e_n)}\sbs I\ot H$. If now
$\sum_{i=1}^ne_ix_i=0$ for some $x_1,\ldots,x_n\in A$, then
$\sum_{i=1}^n\rho(e_i)\rho(x_i)=0$ in $M\ot H$, which shows that
$\rho(x_i)\in\I{\rho(e_1),\ldots,\rho(e_n)}$, and therefore $x_i\in K$ for all
$i=1,\ldots,n$. It follows that $J\sbs K$ by the definition of $J$.
\endproof

\proclaim
Lemma 3.3.
Let $M\in\MAH$ be generated by $e_1,\ldots,e_n$ as an $A$-module. Suppose that
$I$ is an ideal of $A$ such that the ring $A/I\ot H$ is weakly $n$-finite
and $I$ contains no nonzero $H\!$-costable ideals of $A$. If the cosets of
$e_1,\ldots,e_n$ give a basis for the $A/I$-module $M/MI$, then
$e_1,\ldots,e_n$ are a basis for the $A$-module $M$.
\endproclaim

\Proof.
Any relation $\sum_{i=1}^ne_ix_i=0$ implies $x_1,\ldots,x_n\in I$ by the
freeness of $M/MI$. It follows that $\I{e_1,\ldots,e_n}\sbs I$. Applying
Lemma 3.2, we have $K=0$ by the assumptions about $I$. Hence
$\I{e_1,\ldots,e_n}=0$, which means that $e_1,\ldots,e_n$ are linearly
independent over $A$.
\endproof

In the next lemma we denote by $\l V$ the length of $V\in\M_A$ and for each
$Q\in\Max A$ put
$$
r_Q(M)={\l M/MQ\over\l A/Q}.
$$

\proclaim
Lemma 3.4.
Let $A$ be semilocal. Suppose that $M\in\MAH$ is $A$-finite and there exists
$P\in\Max A$ such that $P$ contains no nonzero $H\!$-costable ideals of $A$, the
ring $A/P\ot H$ is weakly finite and $r_P(M)\ge r_Q(M)$ for all $Q\in\Max A$.
Then{\rm:}

\item(i)
a suitable direct sum of finitely many copies of $M$ is a free $A$-module{\rm,}

\item(ii)
if $r_P(M)\in\Z$ then $M$ is a free $A$-module.

\endproclaim

\Proof.
Consider first the case where $r_P(M)=n\in\Z$. By the hypotheses
$$
\l\,(A/Q)^n=n\l A/Q\ge\l M/MQ
$$
for each $Q\in\Max A$. Since $A/Q$ is a simple Artinian ring, we deduce that
the $A/Q$-module $M/MQ$ is an epimorphic image of $(A/Q)^n$, and so $M/MQ$ can
be generated by $n$ elements. Moreover, $M/MP\cong(A/P)^n$ is a free
$A/P$-module. Denote by $J$ the Jacobson radical of $A$. Since
$M/MJ\cong\prod_{Q\in\Max A}M/MQ$, we can find elements $e_1,\ldots,e_n\in M$
whose images generate the $A/Q$-module $M/MQ$ for each $Q$ and give a basis
for the $A/P$-module $M/MP$.  Then $e_1,\ldots,e_n$ generate the $A$-module
$M$ by Nakayama's lemma. Taking $I=P$, we meet the hypotheses of Lemma 3.3
which yields (ii).

In general let $N\in\MAH$ be the direct sum of $t$ copies of $M$ where
$t>0$ is an integer such that $r_P(M)t\in\Z$. We can apply (ii) to $N$.
\endproof

\proclaim
Theorem 3.5.
Suppose that $A$ is a semilocal $H\!$-simple $H\!$-comodule algebra satisfying \C.
Then all objects of $\MAH$ are projective $A$-modules. Moreover, $M\in\MAH$
is a free $A$-module if and only if $M/MQ$ is a free $A/Q$-module for at least
one $Q\in\Max A$.
\endproclaim

\Proof.
If $M$ is $A$-finite then $M$ fulfills the hypotheses of Lemma 3.4 with $P$
taken to be any maximal ideal of $A$ for which $r_P(M)$ attains the maximum
value. Denote by $\V$ the class of all right $A$-modules $V$ such that a
direct sum of finitely many copies of $V$ is a free $A$-module. Lemma 3.4
thus shows that all $A$-finite objects of $\MAH$ belong to $\V$. In
particular, they are direct summands of free $A$-modules, hence projective.
We see also that the family $\F$ of all $\MAH$-subobjects of an arbitrary
$M\in\MAH$ satisfies conditions (a) and (b) of Lemma 2.5 where we take $R=A$.
Indeed, if $N\in\F$ and $N\ne M$, then there exists an $A$-finite $L\in\F$
such that $L\not\sbs N$. Then $N'=N+L$ belongs to $\F$ and properly contains
$N$, while $N'\!/N\in\MAH$ is $A$-finite. Lemmas 2.4 and 2.5 complete the
proof.  In particular, if $M$ is not $A$-finite then $M$ is a free $A$-module.
\endproof

\proclaim
Corollary 3.6.
Under the hypotheses of Theorem {\rm3.5} suppose also that $A/Q$ is a skew
field for some $Q\in\Max A$. Then{\rm:}

\item(i)
all objects of $\MAH$ are free $A$-modules{\rm,}

\item(ii)
$A$ is a simple object of $\MAH${\rm,}

\item(iii)
$A^H$ is a skew field.

\endproclaim

\Proof.
Part (i) follows from the fact that $M/MQ$ is a free $A/Q$-module for any
$M\in\MAH$. If $N\sbs M$ is an $\MAH$-subobject, then $M\cong N\oplus M/N$ in
$\M_A$ where both $N$ and $M/N$ are free $A$-modules. If $0\ne N\ne M$, then
the $A$-module $M$ has a basis consisting of at least two elements. This
proves (ii) since every semilocal ring has invariant base number.
The left multiplication by an element $a\in A$ commutes with the coaction of
$H$ if and only if $a\in A^H$. This shows that $A^H$ is isomorphic to the
endomorphism ring of $A$ as an object of $\MAH$. Now (iii) is a consequence
of (ii) and Schur's Lemma (cf. \cite{Be90, Lemma 2.1}).
\endproof

\proclaim
Proposition 3.7.
Suppose that $A$ is a semilocal $H\!$-comodule algebra satisfying \C and 
having a minimal nonzero $H\!$-costable ideal $M$ which is finitely generated 
in $\M_A$. If at least one maximal ideal of $A$ contains no nonzero 
$H\!$-costable ideals of $A$, then $A$ is $H\!$-simple.
\endproclaim

\Proof.
We may regard $M$ as an $A$-finite object of $\MAH$. If $I$ is any nonzero 
$H\!$-costable ideal of $A$, then $MI\ne0$ because there exists a maximal 
ideal of $A$ containing neither $M$ nor $I$. Then $MI=M$ by the minimality of 
$M$. Denote by $\Om\sbs\Max A$ the subset of those maximal ideals of $A$ which 
contain a nonzero $H\!$-costable ideal of $A$. We see that $MQ=M$ for each 
$Q\in\Om$.  Then $r_Q(M)=0$ whenever $Q\in\Om$. By the hypotheses $\Om\ne\Max 
A$. It follows that the maximum of the numbers $r_Q(M)$, $\,Q\in\Max A$, is 
attained at some maximal ideal $P\notin\Om$. We can now apply Lemma 3.4 which 
shows that $M^t$ is a nonzero free $A$-module for some integer $t>0$. Then 
$(M/MQ)^t$ is a nonzero free $A/Q$-module, so that $MQ\ne M$, for each
$Q\in\Max A$. It follows that $\Om=\varnothing$, and this gives the desired
conclusion.
\endproof

\proclaim
Corollary 3.8.
Let $R$ be a finite dimensional simple algebra, $A$ a finite dimensional
$H\!$-simple $H\!$-comodule algebra. Suppose that $A/P$ is a central simple
algebra for some $P\in\Max A$ and $H$ is weakly finite. Then $R\ot A$ is an
$H\!$-simple $H\!$-comodule algebra with respect to the comodule structure
$\id\ot\rho:R\ot A\to R\ot A\ot H$.
\endproclaim

\Proof.
The algebra $R\ot A/P$ is simple by \cite{Row, Th.~1.7.27}. Hence $P'=R\ot P$
is a maximal ideal of $B=R\ot A$. Take any basis $e_1,\ldots,e_m$ for $R$ and
define $\pi_1,\ldots,\pi_m:B\to A$ such that $x=\sum_{i=1}^me_i\ot\pi_i(x)$
for each $x\in B$. Clearly each $\pi_i$ is a morphism in $\M^H$. Suppose that
$V\sbs P'$ is an $H\!$-subcomodule of $B$. Then $\pi_i(V)\sbs P$ is an
$H\!$-subcomodule of $A$. The ideal of $A$ generated by $\pi_i(V)$ is
$H\!$-costable and is contained in $P$. It follows that $\pi_i(V)=0$ for each
$i=1,\ldots,m$ by the hypotheses. Thus $P'$ contains no nonzero $H\!$-costable
ideals of $B$. Since $\dim B<\infty$, we can apply Proposition 3.7 to the
$H\!$-comodule algebra $B$.
\endproof

\section
4. When are comodule algebras quasi-Frobenius?

Assume throughout this section that $A$ is an $H\!$-comodule algebra.
For finite dimensional algebras we can readily employ duality to derive
additional information. The next lemma generalizes \cite{Ma94, Lemma 2.10}.

\proclaim
Lemma 4.1.
Let $M\in\AMH$. If either $\dim M<\infty$ or $\dim H<\infty$, then
the dual vector space $M^*$ is an object of $\MAH$ in a canonical way.
\endproclaim

\Proof.
Let $\nu:M^*\to\Hom(M,H)$ be the map that assigns to $\xi\in M^*$ the
composite $M\mapr\rho M\ot H\lmapr2{\xi\ot\id}H$. By the hypotheses
$\Hom(M,H)\cong H\ot M^*$, and $\nu$ is a left $H\!$-comodule structure on
$M^*$. Composing $\nu$ with the map $H\ot M^*\to M^*\ot H$ given by
$h\ot\xi\mapsto\xi\ot s(h)$, we obtain a right $H\!$-comodule structure on
$M^*$. If $U\in\M^H$ then the canonical isomorphism
$\Hom(U\ot M,k)\cong\Hom(U,M^*)$ induces a bijective correspondence
between the morphisms $U\ot M\to k$ and $U\to M^*$ in $\M^H$ where $U\ot M$ is
given the tensor product of comodule structures and $k$ is the trivial
$H\!$-comodule.

Clearly $M^*$ is also a right $A$-module with respect to the action
$(\xi a)(m)=\xi(am)$ where $\xi\in M^*$, $\,a\in A$ and $m\in M$. Now there
is a commutative diagram
$$
\diagram{
M^*\ot A\ot M&\lmapr5{\id\ot\la}&M^*\ot M\cr
\diagramskip
\mapd{\mu\ot\id}{}&&\mapd{}{\ev}\cr
\diagramskip
M^*\ot M&\hidewidth\lmapr9{\ev}\hidewidth&k\cr
}
$$
where $\la:A\ot M\to M$ and $\mu:M^*\ot A\to M^*$ are the $A$-module
structures. Both $\la$ and the evaluation map $\ev:M^*\ot M\to k$ are
morphisms in $\M^H$. Hence so too is the composite $\ph=\ev\circ(\id\ot\la)$.
We may regard $\ph$ as a linear map $U\ot M\to k$ with $U=M^*\ot A$. The
diagram shows that $\ph$ corresponds to $\mu:U\to M^*$. It follows that $\mu$
is a morphism in $\M^H$, which yields the required compatibility of module and
comodule structures on $M^*$.
\endproof

\proclaim
Theorem 4.2.
Let $A$ be a finite dimensional $H\!$-simple $H\!$-comodule algebra. If $H$ is
weakly finite then{\rm:}

\item(i)
$A$ is Frobenius,

\item(ii)
all objects of both $\MAH$ and $\AMH$ are projective $A$-modules,

\item(iii)
$M\in\MAH$ {\rm(}resp. $M\in\AMH${\rm)} is a free $A$-module if and only if
$M/MQ$ {\rm(}resp. $M/QM${\rm)} is a free $A/Q$-module for at least one
$Q\in\Max A$.

\endproclaim

\Proof.
By Proposition 2.2 $A$ satisfies \C. Now the assertions concerning $\MAH$ are
restatements of Theorem 3.5. If $M\in\AMH$ is $A$-finite, then $\dim M<\infty$,
and so $M^*$ is a projective right $A$-module in view of Lemma 4.1. Moreover,
$(M^*)^l$ is a free $A$-module for some integer $l\ge0$. In particular, we may
regard $A$ as an object of $\AMH$. Comparing dimensions, we deduce that
$(A^*)^t\cong A^t$ in $\M_A$ for some $t>0$. By the Krull-Schmidt theorem
$A^*\cong A$ in $\M_A$, which yields (i). We have then $A^*\cong A$ in $\AM$
as well, and so the linear duals of free right $A$-modules are free left
$A$-modules. Continuing with an $A$-finite $M\in\AMH$, we conclude that $M^l$
is a free $A$-module for $l$ as above. The proof of (ii) and (iii) is
completed now by Lemmas 2.4 and 2.5.
\endproof

Given $U\in\M_H$ and $V\in\M_A$, we regard $V\ot U$ as a right $A\ot H$-module
with respect to the action $(v\ot u)(a\ot h)=va\ot uh$ where $v\in V$, $u\in
U$, $a\in A$ and $h\in H$. Then $V\ot U$ is a right $A$-module via the
algebra homomorphism $\rho:A\to A\ot H$. Define a right $H\!$-module structure
on the dual of $U$ by the rule $\<u,\xi h\>=\<us(h),\xi\>$ for $u\in U$,
$\xi\in U^*$ and $h\in H$.

\proclaim
Lemma 4.3.
For $U\in\M_H$ and $V,W\in\M_A$ there is a canonical linear injection
$$
\Phi:\Hom_A(W,V\ot U)\to\Hom_A(W\ot U^*,V)
$$
which is a bijection whenever $\dim U<\infty$. More generally, let $X\sbs U$
be any finite dimensional subspace. Then $\Phi$ establishes a bijective
correspondence between the $A$-module homomorphisms $W\to V\ot U$ whose image
is contained in $V\ot X$ and the $A$-module homomorphisms $W\ot U^*\to V$
which factor through $W\ot X^*$.
\endproclaim

\Proof.
Assigning the composite
$W\ot T^*\lmapr3{\th\ot\hbox{}\id}V\ot T\ot T^*\lmapr5{\id\ot\,\ev_T}V$ to
each linear map $\th:W\to V\ot T$, where $T$ is any vector space and
$\ev_T:T\ot T^*\to k$ the evaluation map $t\ot\xi\mapsto\<t,\xi\>$, we obtain
a linear map
$$
\Phi_T:\Hom(W,V\ot T)\to\Hom(W\ot T^*,V).
$$
Clearly $\Phi_T$ is injective. If $\dim T<\infty$, then $\Phi_T$ is bijective.

Note that $\ev_U:U\ot U^*\to k$ is a morphism in $\M_H$, and it follows that
the map $\Phi_U(\th):W\ot U^*\to V$ is a morphism in $\M_A$ whenever so is
$\th:W\to V\ot U$. Since the maps $\Phi_T$ are natural in $T$, we see that
$\Phi_U$ induces a bijection between the linear maps $W\to V\ot U$ whose image
is contained in $V\ot X$ and the linear maps $W\ot U^*\to V$ vanishing on
$W\ot X^\perp$ where $X^\perp$ denotes the kernel of the restriction map
$U^*\to X^*$. Suppose that $\eta:W\ot U^*\to V$ is a morphism in $\M_A$ such
that $\eta(W\ot X^\perp)=0$, and let $\th=\Phi_U^{-1}(\eta)$. One computes
$\th$ by the formula
$$
\th(w)=\sum_{i=1}^n\eta(w\ot e_i^*)\ot e_i\qquad{\rm for\ }w\in W
$$
where $e_1,\ldots,e_n$ is any basis for $X$ and $e_1^*,\ldots,e_n^*\in U^*$
any linear functions whose restrictions to $X$ give the dual basis for $X^*$.
Let us fix $a\in A$ and check that $\th(wa)=\th(w)a$ for every $w\in W$. There
exists a finite dimensional subspace $Z_a\sbs H$ such that $\rho(a)\in A\ot
Z_a$ and $(\rho\ot\id)\circ\rho(a)\in A\ot Z_a\ot Z_a$. Extend the chosen
basis for $X$ to a basis $e_1,\ldots,e_m$ ($m\ge n$) for the subspace
$Y=Xs(Z_a)+X\sbs U$. Lift the elements of the dual basis for $Y^*$ to some
linear functions $e_1^*,\ldots,e_m^*\in U^*$. If $i>n$, then $e_i^*|_X=0$,
whence $\eta(W\ot e_i^*)=0$. If $y\in Y$ then we have
$\sum_{i=1}^m\<y,e_i^*\>e_i=y$. In particular, this holds for elements $y$ in
$X$ and in $Xs(Z_a)$. It follows that
$$
\eqalign{
\sum_{(a)}\sum_{i=1}^mwa_{(0)}\ot\<x,e_i^*a_{(1)}\>\,e_ia_{(2)}
&{}=\sum_{(a)}\sum_{i=1}^mwa_{(0)}\ot\<xs(a_{(1)}\!),e_i^*\>\,e_ia_{(2)}\cr
&{}=\sum_{(a)}wa_{(0)}\ot xs(a_{(1)}\!)a_{(2)}\cr
&{}=wa\ot x=\sum_{i=1}^mwa\ot\<x,e_i^*\>e_i\cr
}
$$
in $W\ot U$ for all $x\in X$. This shows that
$$
\sum_{(a)}\sum_{i=1}^mwa_{(0)}\ot e_i^*a_{(1)}\ot e_ia_{(2)}
-\sum_{i=1}^mwa\ot e_i^*\ot e_i\in W\ot X^\perp\ot U\sbs\Ker\eta\ot U,
$$
and so
$$
\eqalign{
\th(wa)&{}=\sum_{i=1}^n\eta(wa\ot e_i^*)\ot e_i
=\sum_{i=1}^m\eta(wa\ot e_i^*)\ot e_i\cr
&{}=\sum_{(a)}\sum_{i=1}^m\eta(wa_{(0)}\ot e_i^*a_{(1)}\!)\ot e_ia_{(2)}\cr
&{}=\sum_{(a)}\sum_{i=1}^m\eta(w\ot e_i^*)a_{(0)}\ot e_ia_{(1)}\cr
&{}=\sum_{i=1}^m\bigl(\eta(w\ot e_i^*)\ot e_i\bigr)a
=\th(w)a.\cr
}
$$
Thus $\th:W\to V\ot U$ is a morphism in $\M_A$.
\endproof

\proclaim
Lemma 4.4.
Let $U\in\M_H$ and $V\in\M_A$. If $V$ is an injective $A$-module then so too
is $V\ot U$ provided that either $\dim U<\infty$ or $A$ is a right Noetherian
ring, finitely generated as a right $A^H\!$-module.
\endproclaim

\Proof.
We have to show that, whenever $W'$ is a submodule of $W\in\M_A$, every
morphism $\th':W'\to V\ot U$ in $\M_A$ can be extended to a morphism
$\th:W\to V\ot U$. Lemma 4.3 provides $\eta'=\Phi(\th')$ which is a morphism
$W'\ot U^*\to V$ in $\M_A$. As $V$ is injective and $W'\ot U^*$ is an
$A$-submodule of $W\ot U^*$, we can extend $\eta'$ to a morphism
$\eta:W\ot U^*\to V$ in $\M_A$. If $\dim U<\infty$, we have $\eta=\Phi(\th)$
for some morphism $\th:W\to V\ot U$ in $\M_A$, and this $\th$ extends $\th'$.
The same argument works in general as long as we can find an extension $\eta$
vanishing on $W\ot Y^\perp$ for some finite dimensional subspace $Y\sbs U$.

Suppose further that $A$ is a right Noetherian ring, finitely generated as a
right $A^H\!$-module. We may assume that $W'$ is a finitely generated
$A$-module. Indeed, by Baer's criterion \cite{Ka, Th.~5.7.1} it suffices to
consider the case where $W=A$ and $W'$ is a right ideal of $A$. Now $W'$ is
also finitely generated as a right $A^H\!$-module. Let $w_1,\ldots,w_n$ generate
$W'$ over $A^H$. There exists a finite dimensional subspace $X\sbs U$ such
that $\th'(w_i)\in V\ot X$ for all $i=1,\ldots,n$. As
$\th'(w_ia)=\th'(w_i)\cdot(a\ot1)$ for $a\in A^H$, we have $\th'(W')\sbs V\ot
X$. This implies that $\eta'(W'\ot X^\perp)=0$. By a similar argument
$\rho(A)\sbs A\ot Z$ for some finite dimensional subspace $Z\sbs H$. Put
$Y=Xs(Z)\sbs U$. Since $\rho(1)=1\ot1$, we must have $1\in Z$, and therefore
$X\sbs Y$. If $x\in X$, $\xi\in Y^\perp$ and $h\in Z$, then
$\<x,\xi h\>=\<xs(h),\xi\>=0$, which shows that $Y^\perp Z\sbs X^\perp$.
Denote by $N\sbs W\ot U^*$ the $A$-submodule generated by $W\ot Y^\perp$. Then
$$
N\sbs(W\ot Y^\perp)\cdot(A\ot Z)\sbs W\ot X^\perp.
$$
Put $N'=N\cap(W'\ot U^*)$. We see that $N'\sbs W'\ot X^\perp$, and therefore
$\eta'(N')=0$. Now $(W'\ot U^*)/N'$ is embedded into $(W\ot U^*)/N$ as an
$A$-submodule. By injectivity of $V$, the morphism $(W'\ot U^*)/N'\to V$ in
$\M_A$ induced by $\eta'$ extends to a morphism $(W\ot U^*)/N\to V$. The
latter gives the desired $\eta$.
\endproof

\proclaim
Theorem 4.5.
Let $A$ be a semilocal right Noetherian $H\!$-simple $H\!$-comodule algebra
satisfying \C. If $A$ is finitely generated as a right $A^H\!$-module, then
$A$ is a quasi-Frobenius ring, that is, left and right Artinian, left and
right selfinjective.
\endproclaim

\Proof.
For each $V\in\M_A$ we regard $V\ot H$ as an object of $\MAH$ letting $A$
operate via $\rho:A\to A\ot H$ and taking $\id\ot\De:V\ot H\to V\ot H\ot H$ as
a comodule structure. If $M\in\MAH$ then $\rho_M:M\to M\ot H$ is a morphism in
$\MAH$ \cite{Doi84, Example 1.1}. In fact $\rho_M$ is injective as its
composite with $\id\ot\ep:M\ot H\to M$ is the identity transformation of $M$.
Let $E$ be an injective hull of $A$ in $\M_A$. The composite
$$
\ph:A\mapr\rho A\ot H\hookrightarrow E\ot H
$$
is a monomorphism in $\MAH$. Since $(E\ot H)/\ph(A)\in\MAH$ is a
projective $A$-module by Theorem 3.5, $\ph$ is a split monomorphism in $\M_A$.
However, $E\ot H$ is an injective object of $\M_A$ by Lemma 4.4.
Hence so too is $A$. It remains to recall that every right Noetherian right
selfinjective ring is quasi-Frobenius \cite{Ka, Th.~13.2.1}.
\endproof

\section
5. Coactions of finite dimensional Hopf algebras

Let $A$ be an $H\!$-comodule algebra. When $H$ is finite dimensional we can
improve Theorem 4.5 and derive some other conclusions. First comes an
observation valid for any $H$.

\proclaim
Lemma 5.1.
Let $\al:A\to B$ be an algebra homomorphism where $B$ is any associative
algebra, Consider $B\ot H$ as an $H\!$-comodule algebra with respect to
the comodule structure $\id\ot\De:B\ot H\to B\ot H\ot H$. Then\/{\rm:}

\item(i)
$\ph=(\al\ot\id)\circ\rho$ is a homomorphism of $H\!$-comodule algebras
$A\to B\ot H$,

\item(ii)
$\Ker\ph$ is the largest $H\!$-costable ideal of $A$ contained in $\Ker\al$,

\item(iii)
if $\al(A)=B$ then $B\ot H=\ph(A)\cdot(1\ot H)$,

\item(iv)
if $\al(A)=B$ and $H$ has a bijective antipode then
$B\ot H=(1\ot H)\cdot\ph(A)$.

\endproclaim

\Proof.
Clearly $\ph$ is a composite of two algebra homomorphisms. The commutative
diagram
$$
\diagram{
A&\hidewidth\lmapr7\rho\hidewidth&A\ot H&\hidewidth\lmapr{10}
{\al\ot\id}\hidewidth&B\ot H\cr
\diagramskip
\mapd{}\rho&&\mapd{}{\id\ot\De}&&\mapd{}{\id\ot\De}\cr
\diagramskip
A\ot H&\lmapr3{\rho\ot\id}&A\ot H\ot H&\lmapr6{\al\ot\id\ot\id}&B\ot H\ot H\cr
}
$$
shows that $(\ph\ot\id)\circ\rho=(\id\ot\De)\circ\ph$, i.e., $\ph$ respects
the comodule structures as well. Now $\Ker\ph$ is an $H\!$-costable ideal of
$A$. On the other hand, the composite of $\ph$ with $\id\ot\ep:B\ot H\to B$
coincides with $\al$. Hence $\Ker\ph\sbs\Ker\al$. Conversely, $\ph(I)=0$
whenever $I$ is an ideal of $A$ such that $\al(I)=0$ and $\rho(I)\sbs I\ot H$.
Given $a\in A$ and $h\in H$, we have
$$
a\ot h=\sum_{(a)}a_{(0)}\ot a_{(1)}s(a_{(2)}\!)h=\sum_{(a)}\rho(a_{(0)})\cdot
(1\ot s(a_{(1)}\!)h)\in\rho(A)\cdot(1\ot H),
$$\removelastskip
$$
a\ot h=\sum_{(a)}a_{(0)}\ot hs^{-1}(a_{(2)}\!)a_{(1)}=\sum_{(a)}
(1\ot hs^{-1}(a_{(1)}\!))\cdot\rho(a_{(0)}\!)\in(1\ot H)\cdot\rho(A)
$$
in $A\ot H$. Applying the homomorphism $\al\ot\id:A\ot H\to B\ot H$, we deduce
(iii) and (iv).
\endproof

Part (ii) of the next result can be viewed as a generalization of the fact
proved in \cite{Ma92} according to which all Frobenius right coideal
subalgebras of a semisimple Hopf algebra are themselves semisimple.
Corollary 5.3 is close in spirit to Linchenko's result \cite{Li03} which shows
that, under restrictions on $\chr k$, the Jacobson radical of a finite 
dimensional $H$-module algebra for an involutory Hopf algebra is $H$-stable. 
Extension of the latter work to polynomial identity module algebras appeared 
in \cite{Li}. It will be assumed in the rest of this section that
$\dim H<\infty$.

\proclaim
Theorem 5.2.
If $A$ is right Noetherian and there exists $P\in\Max A$ such that
$A/P$ is left Artinian and $P$ contains no nonzero $H\!$-costable ideals of
$A$, then\/{\rm:}

\item(i)
$A$ is $H\!$-simple and is a quasi-Frobenius ring{\rm,}

\item(ii)
$A$ is semisimple Artinian provided that so is $H$.

\endproclaim

\Proof.
Define the algebra homomorphism $\ph:A\to A/P\ot H$ taking $B=A/P$ and
$\al:A\to B$ the canonical homomorphism in Lemma 5.1. Then $\Ker\ph$ is an
$H\!$-costable ideal of $A$ contained in $P$. Hence $\Ker\ph=0$ by the
hypotheses.

We will regard $M=A/P\ot H$ as an $(A/P,A)$-bimodule letting $A$ and $A/P$
operate on $M$ via $\ph$ and the canonical isomorphism $A/P\to A/P\ot1$,
respectively. Since the ring $A/P$ is left Artinian, it has finite length as a
left module over itself (e.g., \cite{Ka, Cor.~9.3.12}). Then so too does the
finitely generated left $A/P$-module $M$. By Lemma 5.1 $M$ is a finitely
generated right $A$-module (as the antipodes of finite dimensional Hopf
algebras are bijective). The latter module is therefore Noetherian. Applying
Lenagan's Theorem \cite{Mc, Th.\ 4.1.6} with left and right sides 
interchanged, we conclude that the right $A$-module $M$ has finite length. 
Since $\ph$ is injective, $A$ has to be right Artinian. Then $A$ is semilocal 
\cite{Ka, Cor.~9.2.3}.  For any ideal $I$ of $A$ the ring $A/I\ot H$ is right 
Artinian as it is finitely generated as a right module over $A/I\ot1$. 
Therefore $A$ satisfies \C. Proposition 3.7 shows that $A$ is $H\!$-simple. To 
complete the proof of (i) we can proceed as in Theorem 4.5 (Lemma 4.4 still 
works).

As in Theorem 4.5 we have $V\ot H\in\MAH$, and so $V\ot H$ is a projective
$A$-module, for any $V\in\M_A$. Suppose that $H$ is semisimple. Then any
cyclic right $H\!$-module $U$ is a direct summand of $H$ in $\M_H$, whence
$V\ot U$ is a direct summand of $V\ot H$ in $\M_A$. Taking $U=k$, we see that
$V\cong V\ot k$ is a direct summand of a free $A$-module. In other words, all
right $A$-modules are projective, whence (ii).
\endproof

\proclaim
Corollary 5.3.
Denote $J=\cap_{P\in\F}P$ where $\F$ is the set of all maximal ideals $P$ of
$A$ such that the factor ring $A/P$ is left Artinian. If $A$ is right
Noetherian and $H$ is semisimple, then $J$ is an $H\!$-costable ideal of $A$.
\endproclaim

\Proof.
For each $P\in\F$ denote by $I_P$ the largest $H\!$-costable ideal of $A$
contained in $P$. The $H\!$-comodule algebra $A/I_P$ has no nonzero
$H\!$-costable ideals contained in $P/I_P$. Theorem 5.2(ii) then shows that
$A/I_P$ is semisimple Artinian. In particular, $I_P$ coincides with the
intersection of those $Q\in\Max A$ for which $I_P\sbs Q$. For each $Q$
appearing here the factor ring $A/Q$ is simple Artinian so that $Q\in\F$.
Hence $J\sbs I_P\sbs P$ for each $P$. It follows that $J=\cap_{P\in\F}I_P$,
and we are done.
\endproof

\proclaim
Proposition 5.4.
Suppose that $A$ is $H\!$-simple and $\dim A<\infty$. Let $V$ be a finite
dimensional and $W$ a simple right $A$-modules. Denoting $D=\End_AW$, we
have
$$
(\dim D)(\dim A)\mid(\dim V)(\dim W)(\dim H).
$$
\endproclaim

\Proof.
Under present hypotheses $D$ is a skew field and $W$ a finite dimensional
vector space over $D$. Let $t=\dim_DW$, and let $Q\in\Max A$ be the
annihilator of $W$ in $A$. Note that $A/Q$ is a simple Artinian ring such that
$A/Q\cong W^t$ in $\M_A$. As before $M=V\ot H$ may be regarded as an object
of $\MAH$. The $A/Q$-module $(M/MQ)^t$ is free since its length is divisible
by $t$. Theorem 3.5 shows that $M^t$ is a free $A$-module, whence $\dim A$
divides $t(\dim M)$. It remains to observe that $(\dim D)t=\dim W$ and $\dim
M=(\dim V)(\dim H)$.
\endproof

Previous result was proved in \cite{Zhu, Th.~2.2} under the additional
assumptions that $k$ is algebraically closed of characteristic zero, $H$ is
semisimple, $A^H=k$ and $V=W$.

\section
6. Finite dimensional coideal subalgebras

Let $A\sbs H$ be a right coideal subalgebra. The opposite multiplication in 
$H$ and the same comultiplication produce a bialgebra $H\op$ which contains 
$A\op$ as a right coideal subalgebra. All conclusions of the next result do 
not change when the pair $A,H$ is replaced with $A\op$, $H\op$. Consequently, 
Theorem 6.1 is valid not only when $H$ is a Hopf algebra but also when $H$ is 
a bialgebra for which $H\op$ has an antipode. Bialgebras with the latter 
property appeared in the literature under the name of {\it anti-Hopf algebras}.
The same extension applies to Theorem 4.2.

\proclaim
Theorem 6.1.
Let $H$ be a weakly finite (anti-)Hopf algebra and $A\sbs H$ a finite 
dimensional right coideal subalgebra. Denote $D=H/HA^+$ and $D'=H/A^+\!H$ 
where $A^+=\Ker\ep|_A$. Then{\rm:}

\item(i)
$A$ is Frobenius and is a simple object of both $\MAH$ and $\AMH$,

\item(ii)
all objects of both $\MAH$ and $\AMH$ are free $A$-modules,

\item(iii)
there are canonical category equivalences $\MAH\approx\M^D$ and
$\AMH\approx\MD$,

\item(iv)
there exist isomorphisms $H\cong D\ot A$ in $\DMA$ and $H\cong A\ot D'$ in
$\ADM$.

\endproclaim

\Proof.
Recall that $H$ is a simple object of $\M_H^H$. If $I$ is an $H\!$-costable
ideal of $A$, then $IH$ is an $\M_H^H$-subobject of $H$, so that $IH$ is
either $0$ or $H$. If, in addition, $I\sbs A^+$ then $\ep(IH)=0$, whence
$IH=0$. We conclude that $A^+$ contains no nonzero $H\!$-costable ideals of $A$.
Proposition 3.7 shows that $A$ is $H\!$-simple. Theorem 4.2 gives the first
part of (i) and also (ii) since $A/A^+\cong k$. Furthermore, $A$ is a simple
object of $\MAH$ by Corollary 3.6. A similar argument works for $\AMH$.

Now $H\in\AMH$ is a free left $A$-module, and the first equivalence in (iii)
is obtained by an application of \cite{Tak79, Th.~1}. If $H$ has a bijective
antipode $s$, then $H\op$ is again a Hopf algebra. As $H\op$ is weakly finite, 
$H\op$ and $A\op$ meet the hypotheses of Theorem 6.1. The second equivalence 
in (iii) is then also fulfilled because $\AMH=\M_{A\op}^{H\op}$ \cite{Doi85, 
(1.2)}. To prove (iii) not assuming the bijectivity of $s$ requires a 
different argument which will be provided in Lemma 6.2. A self-contained proof
of (iv) will be offered in Lemma 6.4; it streamlines the arguments already 
known. Lemma 6.4 can be applied with $M=H$ since $H\in\ADMH$ and $H\in\DMAH$.  
\endproof

\Remarks.
For a finite dimensional $H$ it was shown by Masuoka \cite{Ma92} that (ii) and
(iv) are each equivalent to $A$ being Frobenius, among other equivalent
conditions.

In case of a commutative $H$ the equivalence $\MAH\approx\M^D$ admits an
interpretation in terms of Mackey imprimitivity theory for algebraic groups
and group schemes \cite{Cl83}, \cite{Par80}. The imprimitivity theorem of
Koppinen and Neuvonen \cite{Kop77} for arbitrary finite dimensional Hopf
algebras can be put into this context as well.

Isomorphisms of (iv) were christened the normal basis property. Initially this
property was studied in connection with the structure of Hopf Galois
extensions (e.g, \cite{Kr81}). Some conditions ensuring its fulfillment were
given by Schneider \cite{Sch92}. In particular, \cite{Sch92, Th.~2.4} yields
(iv) in the case where $\dim H<\infty$ and $A$ is a Hopf subalgebra.

According to \cite{MaD92, Prop.~3.2} the first isomorphism in (iv) is
equivalent to the {\it$A$-cocleftness} of $H$, that is, to the existence of a
morphism $H\to A$ in $\M_A$ invertible in the convolution algebra
$\Hom(H,A)$. By \cite{MaD92, Th.~3.4} $H$ is $A$-cocleft provided that $H$
is a faithfully coflat left $D$-comodule and all objects of $\M_{L\ot A}^{L\ot
H}$ are free $L\ot A$-modules where $L$ is an algebraic closure of $k$. The
first of the two hypotheses here means that the cotensor product functor
$?\sq_DH$ is faithfully exact, which is a consequence of the first equivalence
in (iii). Since $L\ot A$ is a finite dimensional right coideal subalgebra of
the Hopf algebra $L\ot H$ over $L$, the freeness in $\M_{L\ot A}^{L\ot H}$ is
also fulfilled by (i). This proves a part of (iv) already at this stage.

Suppose that $\chr k\ne2$. The Hopf algebra described in \cite{Ni89b}, call it
$H$, contains a grouplike element $g$ and a $3\times3$ matrix coalgebra $C$
with a basis $c_{ij}$ ($1\le i,j\le3$) such that
$$
g^2=1,\qquad gc_{ij}=\la_i\la_jc_{ij},\qquad
\De(c_{ij})=\sum_{l=1}^3\,c_{il}\ot c_{lj},\qquad
\ep(c_{ij})=\de_{ij}
$$
where $\la_1=\la_2=1$ and $\la_3=-1$. In this example $C$ is an
$\AMH$-subobject of $H$ but $C$ is not a free left $A$-module. Note that
$s(c_{lt})g=s(gc_{lt})=\la_l\la_ts(c_{lt})$ for $1\le l,t\le3$. Since
$$
\la_i\la_js(c_{lt})c_{ij}=s(c_{lt})gc_{ij}=\la_l\la_ts(c_{lt})c_{ij},
$$
we have $s(c_{lt})c_{ij}=0$ whenever $\la_i\la_j\ne\la_l\la_t$.
Using these equalities together with the identities
$\sum_is(c_{li})c_{ij}=\de_{lj}$ it is easy to see that the matrices
$$
X=\mat{s(c_{11}+c_{13})&s(c_{12})\cr
s(c_{21}+c_{23})&s(c_{22})\cr},
\qquad
Y=\mat{c_{11}+c_{31}&c_{12}+c_{32}\cr
c_{21}&c_{22}\cr},
\qquad
Z=\mat{0&0\cr c_{23}&0\cr}
$$
fulfill the equations $XY=1$ and $XZ=0$ in $\Mat_2(H)$. Thus $X$ has a right
inverse but no left inverse, and so $H$ is not weakly 2-finite.
\endremark

\proclaim
Lemma 6.2.
Let $H$ be any bialgebra, $A$ its right coideal subalgebra and $D=H/HA^+$,
$D'=H/A^+\!H$ quotient coalgebras. Define a functor
$$
\Phi:\MAH\rightsquigarrow\M^D\qquad({\rm resp.,\ }
\Phi:\AMH\rightsquigarrow\MD)
$$
by $M\mapsto M/MA^+$ {\rm(}resp., $M\mapsto M/A^+\!M${\rm)}. If $\Phi$ is
faithfully exact then it is an equivalence. In particular, this is the case
whenever all nonzero objects of $\MAH$ {\rm(}resp., $\AMH${\rm)} are projective
generators in $\M_A$ {\rm(}resp., $\AM${\rm)}.
\endproclaim

The equivalence $\MAH\approx\M^D$ was verified in \cite{Tak79, Th.~1} under
the hypothesis that some $N\in\hbox{}_H\M$ is faithfully flat as a left
$A$-module. If $N$ is such a module then the faithful exactness of $\Phi$
follows from the isomorphisms $M\ot_AN\cong\Phi(M)\ot N$ constructed in
\cite{Tak79} for all $M\in\MAH$. However, the proof of the equivalence given
there does not carry over to our present situation. We will treat $\AMH\!$, 
thereby completing the proof of Theorem 6.1(iii). The other case is similar as
we may replace $A$, $H$ with $A\op$, $H\op$.

\Proof.
The map $\la:H\mapr\De H\ot H\lmapr4{\can\ot\id}D'\ot H$ makes $H$ into a left
$D'$-comodule. Now $H$ is an object of $\AMH$, and $\la$ commutes with left
multiplications by elements of $A$ and the right $H\!$-comodule structure on
$H$. For each vector space $V$ we regard $V\ot H$ as an object of $\AMH$
using available operations on the second tensorand. Define a functor
$\Psi:\MD\rightsquigarrow\AMH$ by the rule
$$
\Psi(V)=V\sqH=\Ker(V\ot H\lmapr8{\id\ot\la-\mu\ot\id}V\ot D'\ot H)
$$
where $\mu:V\to V\ot D'$ is the $D'$-comodule structure on $V\in\MD$ (basic
properties of cotensor products are summarized in \cite{Tak77, Appendix 2};
details can be found in \cite{Da, Ch.~2}). It is immediate that $\Psi(V)$ is
an $\AMH$-subobject of $V\ot H$. For $M\in\AMH$ and $V\in\MD$ there are
natural morphisms
$$
\Xi_M:M\to\Psi\Phi(M),\qquad\Th_V:\Phi\Psi(V)\to V.
$$
Here $\Xi_M$ coincides with the composite
$M\mapr\rho M\ot H\lmapr4{\can\ot\id}\Phi(M)\ot H$; one checks that $\Xi_M$
takes values in $\Phi(M)\sqH$. The map $\id\ot\ep:V\ot H\to V$ vanishes on
$V \ot A^+\!H$, and so the restriction of this map to $\Psi(V)$ factors
through $\Phi\Psi(V)$; one takes $\Th_V$ to be the induced map.

We have $\Psi(D')\cong H$ by \cite{Tak77}. To be precise, $\la$ maps $H$
isomorphically onto $\Psi(D')\sbs D'\ot H$. Since $(\id\ot\ep)\circ\la$
coincides with the canonical projection $H\to D'$, we can identify $\Th_{D'}$
with the identity map $H/A^+\!H\to D'$. In other words, $\Th_{D'}$ is an
isomorphism. It is immediate from the definitions that both $\Phi$ and $\Psi$
commute with arbitrary direct sums. Therefore $\Th_E$ is an isomorphism
whenever $E\in\MD$ is a direct sum of an arbitrary family of copies of $D'$.
Now $D'$ is an injective cogenerator in $\MD$ \cite{Da, Prop.~2.4.3 and
Cor.~2.4.5}. This implies that each $V\in\MD$ is the kernel of a morphism
$E\to E'$ in $\MD$ where both $E$ and $E'$ are direct sums of copies of $D'$.
In the commutative diagram
$$
\diagram{
0&\lmapr3{}\hidewidth&V&\hidewidth\lmapr5{}\hidewidth&E&\hidewidth\lmapr5{}\hidewidth&E'\cr
\diagramskip
&&\mapd{}{\Th_V}&&\mapd{}{\Th_E}&&\mapd{}{\Th_{E'}}\cr
\diagramskip
0&\mapr{}&\Phi\Psi(V)&\mapr{}&\Phi\Psi(E)&\mapr{}&\Phi\Psi(E')\cr
}
$$
both $\Th_E$ and $\Th_{E'}$ are isomorphisms. By \cite{Tak77} the cotensor
products are left exact. In particular, $\Psi$ is left exact. Since $\Phi$ is
exact, the bottom row in the diagram above is exact. It follows that $\Th_V$
is an isomorphism for any $V$.

Let now $M\in\AMH$ and $V=\Phi(M)$. Denote by $K$ and $L$ the kernel and the
cokernel of $\Xi_M$. By the exactness of $\Phi$ we have an exact sequence
$$
0\mapr{}\Phi(K)\mapr{}\Phi(M)\lmapr5{\Phi(\Xi_M)}\Phi\Psi(V)\mapr{}\Phi(L)\mapr{}0.
$$
The composite of $\Phi(\Xi_M)$ with $\Th_V$ is just the identity map
$\Phi(M)\to V$. Since $\Th_V$ is an isomorphism, so is $\Phi(\Xi_M)$ as well.
We deduce that $\Phi(K)=0$ and $\Phi(L)=0$. Then $K=0$ and $L=0$ by
faithfulness of $\Phi$. Thus $\Xi_M$ is an isomorphism for any $M$ as well.

If all objects of $\AMH$ are projective in $\AM$ then all exact sequences in
$\AMH$ split in $\AM$, whence $\Phi$ is exact. If $M\in\AMH$ is a generator in 
$\AM$ then $A^+\!M\ne M$ so that $\Phi(M)\ne0$. This completes the proof.  
\endproof

\proclaim
Lemma 6.3.
Let $D$ be a coalgebra and $M\in\DM$. For each subcoalgebra $C\sbs D$ put
$M_C=\{m\in M\mid\la(m)\in C\ot M\}$ where $\la:M\to D\ot M$ is the given
comodule structure. Then $M\cong D^n$ for some fixed $n\ge0$ if and only if\/
$\dim M_C=n\dim C$ for each finite dimensional subcoalgebra $C$.
\endproclaim

\Proof.
Let $\De:D\to D\ot D$ be the comultiplication and $\ep:D\to k$ the counit.
If $v\in D$ satisfies $\De(v)\in C\ot D$ then $v=(\id\ot\ep)\circ\De(v)\in
C$. Hence $M_C\cong C^n$ whenever $M\cong D^n$. This proves one direction of
the lemma. Conversely, suppose $\dim M_C=n\dim C$ for each finite dimensional
$C$. If $S\sbs D$ is a simple subcoalgebra and $V$ is a simple left
$S$-comodule, then $S\cong V^t$ in $\SM$ where $t=\dim S/\dim V$. It follows
then that $M_S\cong V^{nt}\cong S^n$. The socle of the $D$-comodule $M$
coincides with the sum $\sum M_S$ over all simple subcoalgebras $S\sbs D$.
Hence $\soc M\cong\oplus\,S^n\cong\soc D^n$ in $\DM$. Note that $D^n$ is
an injective $D$-comodule \cite{Da, Cor.~2.4.5}. Any embedding $\soc M\to D^n$
extends therefore to a morphism $\psi:M\to D^n$ in $\DM$. In fact $\psi$ is
injective since so is its restriction to the socle. One has $\psi(M_C)\sbs
C^n$ for each subcoalgebra $C$. Comparing the dimensions, we deduce that
$\psi(M_C)=C^n$ whenever $\dim C<\infty$. As $D$ is a union of finite
dimensional subcoalgebras, $\psi$ is surjective.
\endproof

Let $\DMAH$ be the category whose objects are vector spaces together with a
left $D$-comodule, a right $H\!$-comodule and a right $A$-module structures such
that the $H\!$-comodule and $A$-module structures satisfy the compatibility
condition required in $\MAH$ and these two structures commute with the
$D$-comodule structure. The category $\DAMH$ is defined similarly with
$\AMH$ in place of $\MAH$. For $V\in\DM$ we regard $V\ot A$ as an object of
either $\DMA$ or $\DAM$ so that $A$ operates by multiplications on the second
tensorand and the comodule structure comes from that on $V$.

\proclaim
Lemma 6.4.
Let $A$, $H$ be as in Theorem {\rm6.1} and $D$ any coalgebra. If $M\in\DMAH$ 
{\rm(}resp.\ $M\in\DAMH${\rm)} and $M/MA^+\cong D$ {\rm(}resp.\ 
$M/A^+\!M\cong D${\rm)} in $\DM$ then $M\cong D\ot A$ in $\DMA$ 
{\rm(}resp.\ in $\DAM${\rm)}.
\endproclaim

\Proof.
We consider only $\DMAH$. Since $A$ is Frobenius, there exists $0\ne x\in A$ 
such that $A^+x=0$. For each subcoalgebra $C\sbs D$ define $M_C$ as in Lemma 
6.3. Since the $D$-comodule structure on $M$ commutes with the two other 
structures, $M_C$ is an $\MAH$-subobject of $M$. By Theorem 6.1 $M_C$ is a 
free $A$-module. Then the annihilator of $x$ in $M_C$ coincides with $M_CA^+$, 
and so the action of $x$ induces an isomorphism $M_C/M_CA^+\cong M_Cx$ in 
$\CM$. Taking $C=D$, we obtain an isomorphism $Mx\cong D$ in $\DM$. Next, 
$M_C$ is an $A$-module direct summand of $M$ since $M/M_C\in\MAH$ is a free 
$A$-module. It follows that $M_Cx=M_C\cap Mx\cong C$. If $\dim C<\infty$, 
then $\dim M_C/M_CA^+=\dim C$, and this number is equal to the rank of the 
free $A$-module $M_C$, so that $\dim M_C=(\dim C)(\dim A)$. By Lemma 6.3 
$M\cong D^n$ in $\DM$ where $n=\dim A$.

Suppose that $S\sbs D$ is a simple subcoalgebra and $R=S^*$ the dual simple
algebra. By Corollary 3.8 $B=R\ot A$ is an $H\!$-simple $H\!$-comodule algebra
with a maximal ideal $P=R\ot A^+$. The $S$-comodule structure on objects of
$\SMAH$ corresponds to an $R$-module structure commuting with the $A$-module
and $H\!$-comodule structures. In other words, $\SMAH\cong\M_B^H$.  Now we have
$M_S\in\M_B^H$. As we have seen, $M_S/M_SP=M_S/M_SA^+\cong S$ in
$\SM\approx\M_R$. Every simple finite dimensional algebra is Frobenius. Hence
$S$ is a free $R$-module of rank $1$. By Theorem 4.2 $M_S$ is a free 
$B$-module of rank $1$. There exists then an $R$-submodule, i.e.\ an 
$S$-subcomodule, $U_S\sbs M_S$ such that the map $U_S\ot A\to M_S$ afforded by 
the $A$-module structure is bijective.

We choose such a subcomodule $U_S$ for each simple subcoalgebra $S\sbs D$ and
put $U=\sum_SU_S$. Let $V$ be a maximal $D$-subcomodule of $M$ containing $U$
as an essential subcomodule (so that $U\cap W\ne0$ for every $D$-subcomodule
$0\ne W\sbs V$). Since $M\cong D^n$ is an injective $D$-comodule, $V$ is a
direct summand of $M$ in $\DM$. In particular, $V$ is an injective
$D$-comodule. The canonical map $\ph:V\ot A\to M$ is a morphism in $\DMA$. By the
choice of $U_S$ the restriction of $\ph$ to $U_S\ot A$ is a bijection onto
$M_S$ for each $S$. The sum $\sum_SM_S$ over all simple subcoalgebras is
direct since the $M_S$'s are the isotypic components of the socle $\soc^DM$ of
the $D$-comodule $M$. Hence the restriction of $\ph$ to $U\ot A$ is
injective. As $U\ot A$ is an essential subcomodule of $V\ot A$, we see
that $\ph$ itself is injective. Now $V\ot A$ is a direct sum of copies of $V$
in $\DM$. It is therefore an injective $D$-comodule, whence $\Im\ph$ is a
direct summand of $M$ in $\DM$. On the other hand, the inclusion
$\soc^DM\sbs\Im\ph$ entails the surjectivity of $\ph$. Thus $\ph$ is an
isomorphism. It follows also that $V\cong M/MA^+$ in $\DM$, and we are done.
\endproof

\Remark.
According to \cite{Sch90, Cor.~2.2} (with left and right sides interchanged)
an object $M\in\DMA$ is isomorphic to $D\ot A$ provided that the following two
conditions are fulfilled: $M$ is injective in $\DM$ and $S\sq_DM\cong S\ot A$
in $\SMA$ for each simple subcoalgebra $S\sbs D$. The verification of these
conditions were two main steps in the proof of Theorem 6.4 (note that
$S\sq_DM\cong M_S$).
\endremark

Theorem 6.1 enables us to strengthen \cite{Ma92, Prop.~2.10}:

\proclaim
Corollary 6.5.
If $H$ is a finite dimensional Hopf algebra, then there is a bijective
correspondence between the right coideal subalgebras in $H$ and $H^*$.
\endproclaim

\section
7. Dualization to module algebras

The results of section 3 have their counterparts for $H\!$-module algebras $A$.
In fact condition \C on comodule algebras is no longer needed. The reason is
that the weak finiteness of convolution algebras is recognized by the algebra
argument alone. If $A$ is semilocal then so are all its factor rings $A/I$,
which are therefore weakly finite.

\proclaim
Lemma 7.1.
If $B$ is a weakly finite algebra and $C$ any coalgebra, then the convolution
algebra $\Hom(C,B)$ is weakly finite.
\endproclaim

\Proof.
Note that $\Mat_n\bigl(\Hom(C,B)\bigr)\cong\Hom(C,B_n)$ where $B_n=\Mat_n(B)$.
Two linear functions $C\to B_n$ coincide if and only if they have the same
restriction to every finite dimensional subcoalgebra of $C$. It suffices
therefore to prove that the convolution algebra $\Hom(C,B_n)$ is weakly
1-finite under the assumption that $\dim C<\infty$. In this case $C^*$ is a
finite dimensional algebra and $\Hom(C,B_n)\cong B_n\ot C^*$. The ring $B_n$
is a free module of finite rank over its subring isomorphic to $B$. The same
is valid then for $B_n\ot C^*$, and we may apply Proposition 2.2(e).
\endproof

Let $A$ be an $H\!$-module algebra and $M\in\AH$. Define $\mh\in\Hom(H,M)$ for
each $m\in M$ by the rule $\mh(h)=\ep(h)m$ for $h\in H$. We have then
$(\mh\eta)(h)=m\eta(h)$ for $\eta\in\Hom(H,A)$ and $h\in H$. In Lemmas 7.2,
7.3, 7.4 we assume that {\it$e_1,\ldots,e_n$ generate $M$ as an $A$-module}.

\proclaim
Lemma 7.2.
Under previous assumptions $\eh_1,\ldots,\eh_n$ and
$\tau(e_1),\ldots,\tau(e_n)$ are two systems of generators for the
$\Hom(H,A)$-module $\Hom(H,M)$.
\endproclaim

\Proof.
Let $\xi\in\Hom(H,M)$. There exist $\eta_1,\ldots,\eta_n\in\Hom(H,A)$ such
that $\xi(h)=\sum_{i=1}^ne_i\eta_i(h)$ for all $h\in H$, and it follows that
$\xi=\sum_{i=1}^n\eh_i\eta_i$. Take now $m\in M$ and define $\xi$ by the
rule $\xi(h)=s(h)m$. We get
$$
\eqalign{
\ep(h)m&{}=\sum_{(h)}\,h_{(1)}s(h_{(2)}\!)m
=\sum_{(h)}\,h_{(1)}\bigl(\sum_{i=1}^ne_i\eta_i(h_{(2)}\!)\bigr)\cr
&{}=\sum_{i=1}^n\sum_{(h)}\,(h_{(1)}e_i)\bigl(h_{(2)}\eta_i(h_{(3)}\!)\bigr)
=\sum_{i=1}^n\sum_{(h)}\,(h_{(1)}e_i)\,\th_i(h_{(2)}\!)\cr
}
$$
where $\th_i\in\Hom(H,A)$ is defined by the rule
$\th_i(h)=\sum_{(h)}h_{(1)}\eta_i(h_{(2)}\!)$. This shows that
$\mh=\sum_{i=1}^n\tau(e_i)\,\th_i$, and we have checked already that
$\eh_1,\ldots,\eh_n$ generate $\Hom(H,M)$.
\endproof

\proclaim
Lemma 7.3.
Suppose that $I$ is an ideal of $A$ such that $\I{e_1,\ldots,e_n}\sbs I$ and
the ring $A/I$ is weakly finite. Then there exists an $H\!$-stable ideal $K$ of
$A$ such that $\I{e_1,\ldots,e_n}\sbs K\sbs I$.
\endproclaim

\Proof.
Note that $\Hom(H,I)$ is an ideal of $\Hom(H,A)$, and the factor algebra by
this ideal is isomorphic to $\Hom(H,A/I)$. This factor algebra is weakly
finite by Lemma 7.1. Denote by $K$ the preimage of $\Hom(H,I)$ under
$\tau:A\to\Hom(H,A)$. Clearly $K=\{a\in A\mid Ha\sbs I\}$ so that $K$ is an
$H\!$-stable ideal of $A$ contained in $I$.

If $\eta_1,\ldots,\eta_n\in\Hom(H,A)$ are such that
$\sum_{i=1}^n\eh_i\eta_i=0$, then $\sum_{i=1}^ne_i\eta_i(h)=0$ for all
$h\in H$, and so $\eta_i(H)\sbs I$ for all $i=1,\ldots,n$. It follows that
$\I{\eh_1,\ldots,\eh_n}$ is contained in $\Hom(H,I)$. By Lemma 2.3
$\I{\tau(e_1),\ldots,\tau(e_n)}\sbs\Hom(H,I)$ as well. Suppose that
$\sum_{i=1}^ne_ix_i=0$ for some $x_1,\ldots,x_n\in A$. Then
$\sum_{i=1}^n\tau(e_i)\tau(x_i)=0$ in $\Hom(H,M)$, which shows that
$\tau(x_i)\in\I{\tau(e_1),\ldots,\tau(e_n)}$, yielding the inclusions
$x_i\in K$ for all $i=1,\ldots,n$. This means that $\I{e_1,\ldots,e_n}\sbs K$.
\endproof

The proofs of Lemmas 3.3, 3.4 and Theorem 3.5 generalize without further
complications to the case of module algebras:

\proclaim
Lemma 7.4.
Suppose that $I$ is an ideal of $A$ such that $A/I$ is weakly finite and $I$
contains no nonzero $H\!$-stable ideals of $A$. If the cosets of
$e_1,\ldots,e_n$ give a basis for the $A/I$-module $M/MI$, then
$e_1,\ldots,e_n$ are a basis for the $A$-module $M$.
\endproclaim

\proclaim
Lemma 7.5.
Let $A$ be a semilocal $H\!$-module algebra. Suppose that $M$ is $A$-finite and
there exists $P\in\Max A$ such that $P$ contains no nonzero $H\!$-stable ideals
of $A$ and $r_P(M)\ge r_Q(M)$ for all $Q\in\Max A$ {\rm(}in the notations of
Lemma {\rm3.4)}. Then the conclusions\/ {\rm(i)} and\/ {\rm(ii)} of Lemma
{\rm3.4} are fulfilled.
\endproclaim

\proclaim
Theorem 7.6.
Suppose that $A$ is a semilocal $H\!$-simple $H\!$-module algebra, and let
$M\in\AH$ be locally $A$-finite. Then $M$ is a projective $A$-module.
Moreover, $M$ is a free $A$-module if and only if $M/MQ$ is a free
$A/Q$-module for at least one $Q\in\Max A$.
\endproclaim

The reader may wish to reformulate this result in terms of smash product
algebras. Here $H\cop$ rather than $H$ has to be a Hopf algebra:

\proclaim
Corollary 7.7.
Let $H$ be an anti-Hopf algebra and $A$ a semilocal $H\!$-simple $H\!$-module 
algebra. Then all locally $A$-finite left $A\#H$-modules are projective 
$A$-modules.
\endproclaim

If an $H\!$-module algebra $A$ is not $H\!$-simple, one can ask about projectivity
of localizations. To be precise, let $A$ be right Noetherian, and let $P$ be a
semiprime ideal of $A$. Denote by $\P\sbs A$ the preimage of the set of
regular elements in $A/P$. One says that $P$ is {\it right localizable} if
$\P$ is a right denominator set \cite{Mc, Ch.~4}. Denote by $A_P$ the right
quotient ring of $A$ with respect to $\P$. Suppose that $P$ is right
localizable and contains no nonzero $H\!$-stable ideals of $A$. Is then
$M\ot_AA_P$ a projective $A_P$-module for every locally $A$-finite object
$M\in\AH$? I can prove that this holds true under the assumption that all
prime ideals of $A$ containing $P$ are maximal and $H$ is pointed with
finitely many grouplike elements.

\references
\nextref
Be90
\auth
J.,Bergen;M.,Cohen;D.,Fischman;
\endauth
\paper{Irreducible actions and faithful actions of Hopf algebras}
\journal{Isr. J.~Math.}
\Vol{72}
\Year{1990}
\Pages{5-18}

\nextref
Bou
\auth
N.,Bourbaki;
\endauth
\book{Commutative Algebra}
\publisher{Springer}
\Year{1989}

\nextref
Cl83
\auth
E.,Cline;B.,Parshall;L.,Scott;
\endauth
\paper{A Mackey imprimitivity theory for algebraic groups}
\journal{Math.~Z.}
\Vol{182}
\Year{1983}
\Pages{447-471}

\nextref
Da
\auth
S.,D\u asc\u alescu;C.,N\u ast\u asescu;S.,Raianu;
\endauth
\book{Hopf Algebras, an Introduction}
\bookseries{Pure and Applied Mathematics}
\Vol{235}
\publisher{Marcel Dekker}
\Year{2000}

\nextref
De
\auth
M.,Demazure;P.,Gabriel;
\endauth
\book{Groupes Alg'briques I}
\publisher{Masson}
\Year{1970}

\nextref
Doi84
\auth
Y.,Doi;
\endauth
\paper{Cleft comodule algebras and Hopf modules}
\journal{Comm. Algebra}
\Vol{12}
\Year{1984}
\Pages{1155-1169}

\nextref
Doi85
\auth
Y.,Doi;
\endauth
\paper{Algebras with total integrals}
\journal{Comm. Algebra}
\Vol{13}
\Year{1985}
\Pages{2137-2159}

\nextref
Doi92
\auth
Y.,Doi;
\endauth
\paper{Unifying Hopf modules}
\journal{J. Algebra}
\Vol{153}
\Year{1992}
\Pages{373-385}

\nextref
Dor82
\auth
I.,Doraiswamy;
\endauth
\paper{Projectivity of modules over rings with suitable group action}
\journal{Comm. Algebra}
\Vol{10}
\Year{1982}
\Pages{787-795}

\nextref
Ei
\auth
D.,Eisenbud;
\endauth
\book{Commutative algebra with a view toward algebraic geometry}
\bookseries{Graduate Texts in Math.}
\Vol{150}
\publisher{Springer}
\Year{1995}

\nextref
Ka
\auth
F.,Kasch;
\endauth
\book{Moduln und Ringe}
\publisher{Teubner}
\Year{1977}

\nextref
Kop77
\auth
M.,Koppinen;T.,Neuvonen;
\endauth
\paper{An imprimitivity theorem for Hopf algebras}
\journal{Math. Scand.}
\Vol{41}
\Year{1977}
\Pages{193-198}

\nextref
Kop93
\auth
M.,Koppinen;
\endauth
\paper{Coideal subalgebras in Hopf algebras: Freeness, integrals, smash products}
\journal{Comm. Algebra}
\Vol{21}
\Year{1993}
\Pages{427-444}

\nextref
Kr81
\auth
H.F.,Kreimer;M.,Takeuchi;
\endauth
\paper{Hopf algebras and Galois extensions of an algebra}
\journal{Indiana Univ. Math.~J.}
\Vol{30}
\Year{1981}
\Pages{675-692}

\nextref
Li03
\auth
V.,Linchenko;
\endauth
\paper{Nilpotent subsets of Hopf module algebras}
in ``Groups, Rings, Lie and Hopf Algebras", Kluwer, 2003, pp. 121--127.

\nextref
Li
\auth
V.,Linchenko;S.,Montgomery;L.W.,Small;
\endauth
\paper{Stable Jacobson's radicals and semi\-prime smash products}
\journal{Bull. London Math. Soc.}
\Vol{37}
\Year{2005}
\Pages{860-872}

\nextref
Ma91
\auth
A.,Masuoka;
\endauth
\paper{On Hopf algebras with cocommutative coradicals}
\journal{J. Algebra}
\Vol{144}
\Year{1991}
\Pages{451-466}

\nextref
Ma92
\auth
A.,Masuoka;
\endauth
\paper{Freeness of Hopf algebras over coideal subalgebras}
\journal{Comm. Algebra}
\Vol{20}
\Year{1992}
\Pages{1353-1373}

\nextref
Ma94
\auth
A.,Masuoka;
\endauth
\paper{Coideal subalgebras in finite Hopf algebras}
\journal{J. Algebra}
\Vol{163}
\Year{1994}
\Pages{819-831}

\nextref
Ma94b
\auth
A.,Masuoka;
\endauth
\paper{Quotient theory of Hopf algebras}
in ``Advances in Hopf algebras", Lecture Notes Pure Appl. Math.,
Vol. 158, pp. 107--133, Marcel Dekker, 1994.

\nextref
MaD92
\auth
A.,Masuoka;Y.,Doi;
\endauth
\paper{Generalization of cleft comodule algebras}
\journal{Comm. Algebra}
\Vol{20}
\Year{1992}
\Pages{3703-3721}

\nextref
MaW94
\auth
A.,Masuoka;D.,Wigner;
\endauth
\paper{Faithful flatness of Hopf algebras}
\journal{J. Algebra}
\Vol{170}
\Year{1994}
\Pages{156-164}

\nextref
Mc
\auth
J.C.,McConnell;J.C.,Robson;
\endauth
\book{Noncommutative Noetherian Rings}
\publisher{Wiley}
\Year{1987}

\nextref
Mo83
\auth
S.,Montgomery;
\endauth
\paper{Von Neumann finiteness of tensor products of algebras}
\journal{Comm. Algebra}
\Vol{11}
\Year{1983}
\Pages{595-610}

\nextref
Mo
\auth
S.,Montgomery;
\endauth
\book{Hopf algebras and Their Actions on Rings}
\bookseries{CBMS Regional Conference Series in Mathematics}
\Vol{82}
\publisher{American Mathematical Society}
\Year{1993}

\nextref
Mum
\auth
D.,Mumford;
\endauth
\book{Geometric Invariant Theory}
\bookseries{Ergebnisse der Mathematik und ihrer Grenzgebiete}
\Vol{34}
\publisher{Springer}
\Year{1965}

\nextref
Ni89
\auth
W.D.,Nichols;M.B.,Zoeller;
\endauth
\paper{A Hopf algebra freeness theorem}
\journal{Amer. J. Math.}
\Vol{111}
\Year{1989}
\Pages{381-385}

\nextref
Ni89b
\auth
W.D.,Nichols;M.B.,Zoeller;
\endauth
\paper{Freeness of infinite dimensional Hopf algebras over grouplike subalgebras}
\journal{Comm. Algebra}
\Vol{17}
\Year{1989}
\Pages{413-424}

\nextref
Ni92
\auth
W.D.,Nichols;M.B.,Zoeller;
\endauth
\paper{Freeness of infinite dimensional Hopf algebras}
\journal{Comm. Algebra}
\Vol{20}
\Year{1992}
\Pages{1489-1492}

\nextref
Par80
\auth
B.,Parshall;L.,Scott;
\endauth
\paper{An imprimitivity theorem for algebraic groups}
\journal{Indag. Math.}
\Vol{42}
\Year{1980}
\Pages{39-47}

\nextref
Rad77a
\auth
D.E.,Radford;
\endauth
\paper{Pointed Hopf algebras are free over Hopf subalgebras}
\journal{J. Algebra}
\Vol{45}
\Year{1977}
\Pages{266-273}

\nextref
Rad77b
\auth
D.E.,Radford;
\endauth
\paper{Freeness (projectivity) criteria for Hopf algebras over Hopf subalgebras}
\journal{J. Pure Appl. Algebra}
\Vol{11}
\Year{1977}
\Pages{15-28}

\nextref
Row
\auth
L.H.,Rowen;
\endauth
\book{Ring Theory, Vol. I}
\publisher{Academic Press}
\Year{1988}

\nextref
Sch90
\auth
H.-J.,Schneider;
\endauth
\paper{Principal homogeneous spaces for arbitrary Hopf algebras}
\journal{Isr. J.~Math.}
\Vol{72}
\Year{1990}
\Pages{167-195}

\nextref
Sch92
\auth
H.-J.,Schneider;
\endauth
\paper{Normal basis and transitivity of crossed products for Hopf algebras}
\journal{J. Algebra}
\Vol{152}
\Year{1992}
\Pages{289-312}

\nextref
Sch93
\auth
H.-J.,Schneider;
\endauth
\paper{Some remarks on exact sequences of quantum groups}
\journal{Comm. Algebra}
\Vol{21}
\Year{1993}
\Pages{3337-3357}

\nextref
Sk93
\auth
S.M.,Skryabin;
\endauth
\paper{An algebraic approach to the Lie algebras of Cartan type}
\journal{Comm. Algebra}
\Vol{21}
\Year{1993}
\Pages{1229-1336}

\nextref
Sw
\auth
M.E.,Sweedler;
\endauth
\book{Hopf Algebras}
\publisher{Benjamin}
\Year{1969}

\nextref
Tak72
\auth
M.,Takeuchi;
\endauth
\paper{A correspondence between Hopf ideals and sub-Hopf algebras}
\journal{Manu\-scripta Math.}
\Vol{7}
\Year{1972}
\Pages{251-270}

\nextref
Tak77
\auth
M.,Takeuchi;
\endauth
\paper{Formal schemes over fields}
\journal{Comm. Algebra}
\Vol{5}
\Year{1977}
\Pages{1483-1528}

\nextref
Tak79
\auth
M.,Takeuchi;
\endauth
\paper{Relative Hopf modules---equivalences and freeness criteria}
\journal{J. Algebra}
\Vol{60}
\Year{1979}
\Pages{452-471}

\nextref
Zhu
\auth
Y.,Zhu;
\endauth
\paper{The dimension of irreducible modules for transitive module algebras}
\journal{Comm. Algebra}
\Vol{29}
\Year{2001}
\Pages{2877-2886}

\endreferences
\bye